\documentclass[11pt,a4paper,reqno]{amsart}
\usepackage{graphicx}
\usepackage{amsmath} 
\usepackage{amssymb}
\usepackage{mathtools}
\usepackage{color}
\usepackage{graphicx}
\usepackage{booktabs}
\usepackage{color}
\usepackage{braket}
\usepackage{algorithm}
\usepackage{cases}
\usepackage{array}
\usepackage[labelformat=simple]{subcaption}

\usepackage[shortlabels]{enumitem}
\usepackage[foot]{amsaddr}
\usepackage{amsopn}
\usepackage[nameinlink]{cleveref} 
\usepackage{pgfplots, tikz}
\crefformat{equation}{\textup{#2(#1)#3}}

\newtheorem{theorem}{Theorem}[section]
\newtheorem*{theorem*}{Theorem}
\newtheorem{remark}[theorem]{Remark}
\newtheorem{proposition}[theorem]{Proposition} %
\newtheorem{corollary}[theorem]{Corollary}  %
\newtheorem{definition}[theorem]{Definition}

\crefname{hypothesis}{Hypothesis}{Hypotheses}
\crefname{definition}{Definition}{Definition}
\crefname{subsection}{Subsection}{Subsections}
\crefname{section}{Section}{Sections}
\crefname{theorem}{Theorem}{Theorems}
\crefname{theorem}{Theorem}{Theorems}
\crefname{lemma}{Lemma}{Lemmas}
\crefname{remark}{Remark}{Remark}
\crefname{appendix}{Appendix}{Appendices}
\crefname{proposition}{Proposition}{Propositions}
\crefname{table}{Table}{Tables}
\crefname{figure}{Figure}{Figures}

\newcommand{\dee}{\mathrm{d}}
\newcommand{\OCPT}{(\mathrm{OCP}_T)}
\newcommand{\OCP}{\mathrm{OCP}}
\newcommand{\NLP}{\mathrm{NLP}}
\newcommand{\IVP}{\mathrm{(IVP)}}
\newcommand{\AC}{\mathcal{AC}}
\newcommand{\mbbR}{\mathbb{R}}
\newcommand{\Lcal}{\mathcal{L}}
\newcommand{\zero}{\mathbf{0}}
\newcommand{\U}{\mathcal{U}}
\newcommand{\C}{\mathcal{C}}
\newcommand{\XM}{\mathbb{X}_{\mathcal{M}}}
\newcommand{\dist}{\mathsf{dist}}
\DeclareMathOperator*{\esssup}{ess\,sup}

\newlength{\lowidth}%
\AtBeginDocument{\setlength{\lowidth}{3em}}

\numberwithin{equation}{section}

\allowdisplaybreaks

\title[Nonlinear controlled port-Hamiltonian systems: Existence of (optimal) solutions]{
	Nonlinear controlled port-Hamiltonian systems: Existence of (optimal) solutions
}
\author{Willem Esterhuizen$^{1,2}$}\address{$^1$Optimization-based Control Group, Institute of Mathematics, Technische Universität Ilmenau, Germany\\ Mail: \textsc{\{willem-daniel.esterhuizen, karl.worthmann\}@tu-ilmenau.de}}
\address{$^2$Differential Equations Group, Institute of Mathematics, Technische Universität Ilmenau, Germany}
\author{Bernhard Maschke$^3$}\address{$^3$Universit\'e Claude Bernard Lyon 1, CNRS, LAGEPP UMR 5007, France \\ Mail: \textsc{\{bernhard.maschke\}@univ-lyon1.fr}}
\author{Till Preuster$^{1,4}$}
\author{Manuel Schaller$^{2,4}$}\address{$^4$Faculty of Mathematics, Chemnitz University of Technology, Germany\\ Mail: \textsc{\{till.preuster,manuel.schaller\}@math.tu-chemnitz.de}}
\author{Karl Worthmann$^1$}

\thanks{W.E., K.W. and M.S.\ gratefully acknowledge funding by the Deutsche Forschungsgemeinschaft (DFG, German Research Foundation) -- Project-ID 519323897.}
\thanks{This is the authors' version of the article published in SIAM Journal on Control and Optimization, vol.~64, no.~4, pp.~2569--2592, 2026. DOI: https://doi.org/10.1137/24M1705093. \copyright~2026 Society for Industrial and Applied Mathematics.}

\begin{document}

\begin{abstract}
	We investigate the existence of solutions of reversible and irreversible port-Hamilto\-nian systems. To this end, we utilize the associated exergy, a function that is composed of the system's Hamiltonian and entropy, to prove global existence in time for bounded control functions. 
	Then, we rigorously verify our existence conditions for two examples, the gas-piston system and a network of heat exchangers.
	Last, we explore model predictive control tailored to irreversible port-Hamiltonian systems by means of a numerical case study with a heat exchanger network.
\end{abstract}

\maketitle


\smallskip
\noindent \textbf{Mathematics subject classications.}  
34A12, 
49J15, 
80M50, 
93D20 

\bigskip

\section{Introduction}
	\noindent Port-Hamiltonian systems (PHSs) 
	provide a modular modeling framework that is especially convenient for the modelling, \cite{hauschild2020port, lohmayer2021exergetic, warsewa2021port}, analysis \cite{jacob2012linear} and control \cite{van2014port, ramirez2016passivity, ortega2001putting}, of multi-physical systems. At the heart of port-Hamiltonian modeling is the representation of the structure of physical systems, such as its topology, for instance the Kirchhoff's laws or the energy of the system and the associated co-energy variables and the dissipative relations.
	
	While feedback control leveraging the port-Hamiltonian structure and passivity is a well-developed topic~\cite{van2000l2,van2014port}, there are significantly less works considering optimal control of port-Hamiltonian systems. 
	In \cite{kolsch2021optimal}, the authors present an approximate dynamic programming approach (where an optimal control problem's value function is learnt during online operation for use in specifying a feedback law) that uses the system's Hamiltonian as an initial guess. 
	The papers~\cite{sprangers2014reinforcement} and~\cite{nageshrao2014passivity} show how the port-Hamiltonian structure can be exploited in reinforcement learning. 
	The work \cite{sato2017riemannian} considers linear PHSs and presents a structure-preserving $\mathcal{H}_2$-controller design approach in the sense that the closed-loop system is again port-Hamiltonian. 
	Similarly, in \cite{breiten2023structure} a structure-preserving $\mathcal{H}_\infty$ design is studied.
	In \cite{wu2020reduced} the authors present an approach to define a finite-dimensional linear-quadratic-Gaussian (LQG) controller via a structure-preserving reduced-order model of an infinite-dimensional PHS. 
	
	Moreover, many recent papers consider the open-loop optimal control of PHSs. 
	This includes the contribution \cite{schaller2021control} where the authors focus on linear dissipative port-Hamiltonian systems (DPHSs) and characterise the solution to the problem of minimising the supplied energy. One of the paper's main results states that the solution, which may contain \emph{singular arcs}, exhibits \emph{turnpike behaviour}, with the optimal state staying close to a \emph{non-dissipative subspace} for most of the time over the problem's finite horizon. 
	The paper \cite{faulwasser2023hidden} presents an in-depth study of this problem's singular solutions, and a relation to general singular optimal control and the analysis via Goh condition was given in \cite{soledad2024conditions}. 
	The paper \cite{faulwasser2022optimal} extends the study to descriptor systems (those with additional algebraic constraints). 
	In~\cite{karsai2024manifold} the author further studies the turnpike phenomenon of these energy-minimising problems, focussing on nonlinear port-Hamiltonian descriptor systems. 
    
	In this paper, we consider well-posedness of the class of reversible-irreversible port-Hamiltonian systems \cite{ramirez2022overview}, suited to represent irreversible thermodynamic systems such as chemical reactors or heat exchangers. These nonlinear systems are a~generalization of dissipative port-Hamiltonian systems and are generated by two functions: The Hamiltonian function corresponds to the total energy, while its counterpart - the entropy functional - represents the total entropy of the system, see, e.g.,~\cite{ramirez2013modelling}. In~\cite{philipp2023optimal}, an optimization problem modeling the minimisation of energy supply, the entropy growth and the exergy production is formulated and the authors analyze
	turnpike behaviour with respect to the set of thermodynamic equilibria. Note that metriplectic and GENERIC systems~\cite{ottinger2005beyond} may be transformed into reversible-irreversible port-Hamiltonian systems~\cite{Goreac_ArXiV_2024}, hence broadening the scope of our results. 
	
	While there are various works analyzing the qualitative behavior of optimal solutions, 
	the topic of existence of solutions to PH optimal control problems (OCPs) has been largely unexplored. 
	An exception is the paper \cite{doganay2023optimal} where the authors consider DPHSs with quadratic Hamiltonian and present an existence result for a finite-horizon OCP with running and terminal cost that are both continuous and convex. 
	
	The contributions of this paper are as follows. 
	In Theorem~\ref{thm:IVP_main}, we present conditions under which a global solution to nonlinear reversible-irreversible port-Hamiltonian initial-value problem exists.
	In our reasoning, we leverage the exergy function together with the energy and entropy balance inherent to this class of systems
	to show that the unique local solution remains within a compact set over its maximal interval of existence, thus implying existence on the whole time axis.
	As a consequence, as dissipative Hamiltonian systems are a particular case of the considered class, we provide a global existence result for nonlinear dissipative port-Hamiltonian systems in Corollary~\ref{cor:dissipative}.
	Using this existence result we derive an existence result to nonlinear energy- and entropy-optimal control problems in Theorem~\ref{thm:OCP_existence}.

	\textbf{Outline.} 
	 In Section~\ref{sec:pH_systems} we introduce the model class of reversible-irreversible port-Hamiltonian systems (RI-PHSs). 	Section~\ref{sec:IVP} is concerned with the existence of solutions to initial-value problems for RI-PHSs and presents the paper's first contribution (Theorem~\ref{thm:IVP_main}) whereas Section~\ref{sec:OCP_existence} considers existence questions for OCPs and presents the paper's second contribution (Theorem~\ref{thm:OCP_existence}).  Then, we apply the abstract existence results for two specific examples of RI-PHS: heat exchanger networks (Section~\ref{sec:app:heat_ex}) and a gas-piston systems (Section~\ref{sec:app:gaspiston}). 	In Section~\ref{sec:numerics} we provide an extensive numerical case study involving a network of heat exchangers with particular focus on stability of solutions and provide an outlook towards entropy- and energy-oriented Model Predictive Control. 
	We then briefly conclude the paper in Section~\ref{sec:conclusion}.

	\textbf{Notation.} Throughout the paper, 
	we denote the finite-dimensional $p$-norm by $\| \boldsymbol\cdot \|_p$, $p\in[1, \infty) \cup \{\infty\}$. If clear from context or if the particular norm is irrelevant due to equivalence, we simply write  $\| \boldsymbol\cdot \|$.
	Given $I\subseteq \mbbR$, an arbitrary interval, we let $\mathcal{L}(I,\mbbR^m)$ denote the space of Lebesgue-measurable functions mapping $I$ to $\mbbR^m$. 
	Given two functions $u:I\rightarrow \mbbR^m$ and $y:I\rightarrow \mbbR^m$, set
	$\| u \|^p_{\Lcal^p} \coloneqq \left(\int_I \| u(s)\|^p\,\, \dee s\right), 	\| u \|_{\Lcal^{\infty}} \coloneqq \esssup_{t\in I}\| u(t)\|, \langle u, y\rangle \coloneqq \int_I  u(s)^\top y(s)\,\, \dee s.$ 
	We let $\mathcal{L}^p(I,\mbbR^m)$, $p\in[1, \infty) \cup \{\infty\}$, denote all $u\in\mathcal{L}(I, \mbbR^m)$ for which $\| u \|_{\Lcal_p} < \infty$, and let $\Lcal^{1}_{\mathrm{loc}}(I,\mbbR^m)$ denote all $u\in\mathcal{L}(I, \mbbR^m)$ which are locally Lebesgue integrable, that is, integrable on every compact subset of $I$. 
	With $I\subset\mbbR$ a compact interval, we let $\AC(I,\mbbR^n)$ denote all absolutely continuous functions mapping $I$ to $\mbbR^n$. 
	If $I\subseteq\mbbR$ is an interval, not necessarily compact, a function is said to be locally absolutely continuous if it is absolutely continuous on every compact subset of $I$. 
	With $N\in\mathbb{N}$, we define $[1:N] := [1,N] \cap \mathbb{N}$, 
	that is, the set of integers from 1 to $N$. 
	Given a differentiable function $v:\mbbR^n\rightarrow \mbbR$, its gradient is denoted by $v_x:\mbbR^n\rightarrow \mbbR^n$ and for a twice differentiable function $v:\mbbR^n\rightarrow \mbbR$, the Hessian is denoted by $v_{xx}:\mbbR^n\rightarrow \mbbR^{n \times n}$. We indicate (twice) differentiable functions by $C^1(\mbbR^n, \mbbR),C^2(\mbbR^n, \mbbR)$, respectively.
	\sloppy
	Given two differentiable functions $v,w:\mbbR^n\rightarrow \mbbR$ and a skew-symmetric matrix $J\in\mbbR^{n\times n}$, the \emph{Poisson bracket} is defined as 
	$\{v,w\}_J(x) := v_x(x)^\top J w_x(x),\quad x\in\mbbR^n.$
	\fussy
	A positive definite matrix $Q\in\mbbR^{n\times n}$ is indicated by $Q\succ 0$. 
	Its largest (smallest) eigenvalue is denoted by $\overline\sigma(Q)$ ($\underline \sigma(Q)$). 

\section{Coupled Reversible-Irreversible Port-Hamiltonian Systems}\label{sec:pH_systems}
	\noindent In this section, we introduce the considered class of nonlinear control systems as proposed in~\cite{ramirez2013modelling} and we provide two applications that will serve as running examples in this work. For a broad overview over the topic of (reversible-)irreversible port-Hamiltonian systems, we refer to the survey article~\cite{ramirez2022overview}.	

\sloppy
    
	\begin{definition}\label{def:RIPHS}
		A reversible-irreversible port-Hamiltonian system (RI-PHS) is defined by\\
        \hspace*{1em}(i) a matrix-valued locally Lipschitz continuous function $J_0 : \mbbR^n \rightarrow \mbbR^{n\times n}$, which is \\\hspace*{2em}pointwise skew-symmetric, that is, $J_0(x) = -J_0(x)^\top$ for all $x\in \mbbR^n$,\\[.5mm]
		\hspace*{1em}(ii) $N$ constant skew-symmetric matrices $(J_k)_{k \in [1:N]} \in\mbbR^{n\times n}$,\\[.5mm]
		\hspace*{1em}(iii) a function $H\in C^1(\mbbR^n, \mbbR)$ with locally Lipschitz gradient $H_x : \mbbR^n \rightarrow \mbbR^n$,\\[.5mm]
		\hspace*{1em}(iv) a function $S\in C^1(\mbbR^n, \mbbR)$ with locally Lipschitz gradient $S_x : \mbbR^n \rightarrow \mbbR^n$, which 
        is a \hspace*{2em}Casimir function of $J_0$, that is, $J_0(x) S_x(x) = 0$ for all $x\in\mbbR^n$,\\[.5mm]
		\hspace*{1em}(v) $N$ positive and locally Lipschitz functions $(\gamma_k)_{k \in [1:N]}:\mbbR^n \times \mbbR^n \rightarrow \mbbR_{> 0}$,\\[.5mm]
		\hspace*{1em}(vi) a matrix-valued locally Lipschitz function $g:\mbbR^n \times \mbbR^n \rightarrow \mbbR^{n \times m}$,\\[.5mm]
		\hspace*{1em}(vii) the state equation,
			\begin{align}\label{eq:RIPH}
			\dot x & =  \left( J_0(x) \!+\! \sum\nolimits_{k=1}^N \gamma_k(x, H_x(x)) \{S, H\}_{J_k}(x) J_k \right)\! H_x(x) + g(x, H_x(x)) u, 
			\end{align}
             where $u \in \U([0,\infty), \mbbR^m)$ denotes a suitable control (we clarify the notation $\U([0,\infty), \mbbR^m)$ in Section~\ref{sec:IVP}.) \\[.5mm]
		\hspace*{1em}(viii) energy- resp.\ entropy-conjugated outputs $y_H, y_S\in\mbbR^m$, given by,
			\begin{align}\label{eq:def_ys}
			y_H := g(x, H_x)^\top H_x,\quad y_S := g(x, H_x)^\top S_x.
			\end{align}
	\end{definition}
In the preceding definition the introduced functions $H$ and $S$ are called the \textit{Hamiltonian} and \textit{entropy functional}, respectively.
	Straightforward computations using skew-symmetry of the interconnection matrices $(J_k)_{k \in [1:N]}$ show that any solution to \eqref{eq:RIPH} (the precise notion of solutions will be defined in Section~\ref{sec:IVP}) obeys the \emph{power balance} 
	\begin{equation}
	\frac{\dee H}{\dee t}(x(t)) =  H_x(x(t))^\top \mathcal{J}(x(t)) H_x(x(t)) +  u(t)^\top  y_H(t) = u(t)^\top  y_H(t) , \label{eq:energy_balance}
	\end{equation}
	as $\mathcal{J}(x) := \left( J_0(x) + \sum\nolimits_{k=1}^N \gamma_k(x, H_x(x)) \{S, H\}_{J_k}(x) J_k \right)$  is skew-symmetric, and the \emph{entropy balance} 
	\begin{align}
	 \label{eq:entropy_balance}
	\begin{split}
	\frac{\dee S}{\dee t}(x(t)) &= \sum\nolimits_{k=1}^N \!\gamma_k(x(t), H_x(x(t))) \left(\{S, H\}_{J_k}(x(t))\right)^2  \!+\! u(t)^\top y_S(t) \geq  u(t)^\top y_S(t),
	\end{split}
	\end{align}
	which follows from the positivity of the $(\gamma_k)_{k \in [1:N]}$. 
	For a closed system, i.e., $u \equiv 0$, \eqref{eq:energy_balance} shows that the first law of thermodynamics, the conservation of energy, is satisfied. Correspondingly, Equation \eqref{eq:entropy_balance} expresses the second law of thermodynamics, that is, the \emph{irreversible} entropy growth.  We will provide examples of two central applications, i.e., a heat exchanger network and a gas-piston system, in Section~\ref{sec:app:heat_ex} in Section~\ref{sec:app:gaspiston} respectively.

	\section{{Reversible-irreversible pH Systems}: Existence of Solutions} 
    \label{sec:IVP}
	\noindent In this part, we provide our first main result ensuring global existence of solutions to the initial-value problem (IVP) 
	\begin{numcases}{\IVP:\quad}
	\dot x(t) & $ = \quad f(x(t),u(t))\quad \mathrm{a.e.}\quad t\in[0,\infty)$ \nonumber \\
	x(0) & $= \quad x_0$, \nonumber \\
	u   & $\in \quad\U([0,\infty),\mbbR^m),$ \nonumber
	\end{numcases}
	where $f:\mbbR^n\times \mbbR^m \rightarrow \mbbR^n$ has the particular structure as in Definition \ref{def:RIPHS}, that is
	\begin{align*}
	f(x,u) = \bigg( J_0(x) + \sum\nolimits_{k=1}^N \gamma_k(x, H_x(x)) \{S, H\}_{J_k}(x) J_k \bigg) H_x(x) + g(x, H_x(x)) u,
	\end{align*}
	and where $x_0\in\mbbR^n$ is the initial state. In the following, $\U([0,\infty), \mbbR^m)$ serves as a placeholder for the space of control functions, that is,
	\begin{align*}
	\U([0,\infty), \mbbR^m)  := \Lcal^{\infty}([0,\infty),\mbbR^m)  \quad \mathrm{or} \quad \U([0,\infty), \mbbR^m):= \Lcal^{1}([0,\infty),\mbbR^m). 
	\end{align*}
	Following \cite[App.\ C]{Sontag2013}, we adopt the usual notion of solutions to $\IVP$ and  provide a well-known (local) existence result.
	\begin{proposition}\label{prop:local_sol}
		Given $x_0\in\mbbR^n$ and $u\in \U([0,\infty), \mbbR^m)$. Then, there exists a unique solution to $\IVP$ on the maximal interval of existence $I_{x_0,u} = [0,t_{x_0,u}) \subseteq [0,\infty)$, with $t_{x_0,u} > 0$, that is, there is a unique  locally absolutely continuous function $x:J \rightarrow \mbbR^n$ that satisfies the integral equation
		\[
		x(t) = x_0 + \int_0^t f(x(s),u(s))\,\, \dee s \qquad \forall t\in J.
		\]
	\end{proposition}
	\begin{proof}
		The function $f$ is locally Lipschitz in $x$. Further, as $\Lcal^{\infty}([0,\infty),\mbbR^m) \subset \Lcal^{1}_{\mathrm{loc}}([0,\infty),\mbbR^m)$ and $\Lcal^{1}([0,\infty),\mbbR^m) \subset \Lcal^{1}_{\mathrm{loc}}([0,\infty),\mbbR^m)$, \cite[Thm~54, App.~C]{Sontag2013} may be applied for every $x_0\in\mbbR^n$ and $u\in\U([0,\infty),\mbbR^m)$, which yields the claim.
	\end{proof} 	
	
	\subsection{A global existence result for RI-PHS}
	We now present the paper's first contribution, providing conditions under which a \emph{global} solution to $\IVP$ exists, that is, one for which $I_{x_0, u} = [0,\infty)$.
	To that end we will consider the \emph{exergy function}
	\begin{equation}
	E(x) := H(x) - T_0S(x), \label{def:exergy}
	\end{equation}
	where $T_0 >0$ is a reference temperature (so that $E$ and $T_0S(x)$ have the units of energy as $H$ has). 
	The essence of the proof is that under suitable conditions on $E$, its gradient $E_x$ and the input matrix $g$, the maximal solution is constrained to sublevel sets of $E$ that are compact. 
	
    The underlying idea of using the exergy function as a Lyapunov-type functional was already leveraged before in numerous applications, e.g.\ for heat conduction in rigid bodies~\cite{bulivcek2019thermodynamics} and flow problems~\cite{skiba2024exergy} see also~\cite{sangi2019application} for a survey focusing on control design. Regarding existence of solutions for nonlinear partial differential equations in quasistatic electromagnetism, energy estimates and the port-Hamiltonian structure were central in~\cite{chill2023analysis}. In this work, we leverage the exergy to show existence of solutions for a class of port-Hamiltonian systems with particular focus on optimal control.

    We briefly outline the logic of the existence proof, which is conducted by contradiction. To this end, assume that the maximal interval of existence $I_{x_0,u} = [0,t_{x_0,u})$ fulfills $t_{x_0,u}< \infty$. We show that under suitable assumptions on the exergy, the solution remains bounded on $I_{x_0,u}$. By the continuation result in \cite[Prop~C.3.6]{Sontag2013}, this implies that the solution can be extended beyond $t_{x_0,u}$. This would contradict maximality and hence, $t_{x_0,u}=\infty$.

		\begin{theorem}[Global existence of solutions]\label{thm:IVP_main}
		Consider $\IVP$ and suppose the following assumptions hold:
		\noindent
		\begin{enumerate}
			\item[(A1)] The exergy function $E$ defined in \eqref{def:exergy} is unbounded along unbounded local solutions, that is,
			\[
			\lim_{t\rightarrow t_{x_0,u}} \left| \left| x(t)  \right|\right|  = \infty \qquad \Longrightarrow \qquad  \lim_{t\rightarrow t_{x_0,u}} E(x(t)) = \infty. 
			\]
			\item[(A2)]
			Either:
			\begin{enumerate}
				\item[(A2.1)] $\U([0,\infty),\mbbR^m) = \Lcal^{1}([0,\infty),\mbbR^m)$ and there exist a $p\in(1,\infty)$ and a constant $C > 0$ such that for all $t \in I_{x_0,u}$:
				\[
				\left| \left| g(x(t), H_x(x(t)))^\top E_x(x(t)) \right|\right|^p  \leq C E(x(t)).
				\]
			\end{enumerate}
			Or:
			\begin{enumerate}
				\item[(A2.2)] $\U([0,\infty),\mbbR^m) = \Lcal^{\infty}([0,\infty),\mbbR^m)$ and there exists a constant $C > 0$ such that for all $t \in I_{x_0,u}$:

				\[
				\left| \left| g(x(t), H_x(x(t)))^\top E_x(x(t)) \right|\right|  \leq C E(x(t)).
				\]
			\end{enumerate}
		\end{enumerate}
		Then $I_{x_0,u} = [0,\infty)$. 
	\end{theorem}
	\begin{proof}
        From the energy and entropy balance, see equations \eqref{eq:energy_balance} and \eqref{eq:entropy_balance} we have
		\begin{align*}
		\frac{\dee E}{\dee t}(x(t)) = \frac{\dee H}{\dee t}(x(t)) \!-\! T_0\frac{\dee S}{\dee t}(x(t)) 
		& \leq \!u(t)^\top\! g(x(t), H_x(x(t)) )^\top \!\left( H_x(x(t)) \!- \!T_0S_x(x(t)) \right) \\&=  \!u(t)^\top\! g(x(t), H_x(x(t)) )^\top E_x(x(t)).
		\end{align*}
		Assume (A2.1) holds. 
		Then, using Hölder's inequality, we get for any $t\in I_{x_0,u}$
		\begin{align*}
		E(x(t)) &\leq E(x_0) +  \int\nolimits_0^t  u(s)^\top g(x(s), H_x(x(s)))^\top E_x(x(s))\,\,\dee s  \nonumber \\
		& \leq E(x_0) \!+\!| \langle u, g(x, H_x(x))^\top E_x(x) \rangle | \nonumber  \leq E(x_0) \!+ \!\| u \|_{\Lcal^{1}} \| g(x, H_x(x))^\top \!E_x(x)\|_{\Lcal^{\infty}}. 
		\end{align*}
		Let $q\in(1,\infty)$ be a number satisfying $\frac{1}{q} + \frac{1}{p} = 1$ where $p$ is such that (A2.1) holds. 
		Then from Young's inequality we obtain
		\begin{align}
		E(x(t)) &\leq E(x_0) + \frac{\| u \|_{\Lcal^{1}}^q}{q\varepsilon^q} + \frac{\varepsilon^p\| g^\top(x, H_x(x))E_x(x) \|_{\Lcal^{\infty}}^p}{p} \label{ineq:E_after_young},
		\end{align}
		where $\varepsilon > 0$ is an arbitrary number. 
		 Focussing on the third term in \eqref{ineq:E_after_young} we observe, by (A2.1), that
		\begin{align*}
		\| g^\top(x, H_x(x))E_x(x) \|_{\Lcal^{\infty}}^p  &= \left[ \sup_{s \in I_{x_0,u} }\|g^\top(x(s), H_x(x(s)))E_x(x(s)) \| \right]^p \\ 
        &=  \sup_{s\in I_{x_0,u} }\left[\|g^\top(x(s), H_x(x(s)))E_x(x(s)) \|^p \right] \\
        &\leq \sup_{s\in I_{x_0,u}} C E(x(s)). 
		\end{align*}
        We now substitute this upper approximation on $\| g^\top(x, H_x(x))E_x(x) \|_{\Lcal^{\infty}}^p$ into \eqref{ineq:E_after_young} to obtain,
		\begin{align}
		E(x(t)) &\leq E(x_0) + \frac{\| u \|_{\Lcal^{1}}^q}{q\varepsilon^q} + \frac{\varepsilon^p}{p}C \sup_{s\in I_{x_0,u}}E(x(s)). \label{ineq:E_after_young_and_A21}
		\end{align}
		Now we take the supremum over $t$ on the left-hand side of \eqref{ineq:E_after_young_and_A21} (note that the right-hand side is independent of $t$) and hence
		\begin{align*}
		\sup_{s\in I_{x_0,u}}E(x(s)) &\leq E(x_0) + \frac{\| u \|_{\Lcal^{1}}^q}{q\varepsilon^q} + \frac{\varepsilon^p}{p}C \sup_{s\in I_{x_0,u}}E(x(s)).
		\end{align*}
		Let $\varepsilon$ satisfy $0 <\varepsilon < \left(\frac{p}{C}\right)^{1/p}$. 
		Thus, $1 - \frac{\varepsilon^p}{p}C > 0$ and
		\[
		E(x(t)) \leq \frac{E(x_0) +\frac{\| u \|_{\Lcal^{1}}^q}{q\varepsilon^q}}{1 - \frac{\varepsilon^p}{p}C} \eqqcolon  \bar C
		\]
		for all $ t\in I_{x_0,u}$.
		We note that $\bar C < \infty$ because $u\in\Lcal^1([0,\infty),\mbbR^m)$.
		The contrapositive of the statement in (A1) states that for any $\bar C < \infty$ there is $C_2(\bar C) < \infty$ such that
		\[
		\lim_{t\rightarrow t_{x_0,u}}  E(x(t)) \leq \bar C \qquad \Longrightarrow \qquad \lim_{t\rightarrow t_{x_0,u}}  \| x(t) \| \leq C_2(\bar C).
		\]
		Therefore, for all $t\in I_{x_0,u}$,
		\[
		\| x(t) \| \leq C_2(\bar C) < \infty 
		\]
		meaning that $x(t)$ is constrained to a compact set over $I_{x_0, u}$ and we can conclude by \cite[Prop~C.3.6]{Sontag2013} that $I_{x_0,u} = [0,\infty)$. 
		
		Assume now that (A2.2) holds instead. 
		Again using Hölder's inequality, we get for any $t\in I_{x_0, u}$, 		\begin{align*}
		E(x(t)) &  \leq E(x_0) + \| u \|_{\Lcal^{\infty}} \| g^\top(x, H_x(x))E_x(x)\|_{\Lcal^1} \nonumber \\
		& = E(x_0) + \| u \|_{\Lcal^{\infty}} \int_0^t\| g^\top(x(s), H_x(x(s))) E_x(x(s))\|\,\, \dee s \nonumber \\
		& \leq E(x_0) + \| u \|_{\Lcal^{\infty}} C \int_0^t E(x(s)) \,\, \dee s.  \nonumber
		\end{align*}
        Since we assume $t_{x_0,u}<\infty$, Grönwall's inequality \cite[Section III.1]{hartman2002ordinary}, yields
        \begin{equation*}
            E(x(t)) \leq E(x_0)e^{C \| u \|_{\Lcal^\infty}t} \eqqcolon \bar C(t) \leq \bar C(t_{x_0,u})
        \end{equation*}
        for any $t\in I_{x_0, u}$. With that we observe
        \begin{equation*}
            \sup_{s \in I_{x_0,u}} E(x(s)) \leq \bar C(t_{x_0,u}) < \infty,
        \end{equation*}
        where the boundedness of  $\bar C(t_{x_0,u})$ is a consequence of $u\in\Lcal^{\infty}([0,\infty),\mbbR^m)$. Suppose by contradiction that $\sup_{s \in I_{x_0,u}} \|x(s)\|=\infty$. Then there exists a sequence $t_n \nearrow t_{x_0,u}$ such that $\|x(t_n)\| \to \infty$. Consequently, (A1) implies $E(x(t_n)) \to \infty$ contradicting the boundedness of $E$ on $I_{x_0,u}$. Therefore, for any finite $t_{x_0,u} < \infty$, the trajectory $x(t)$ remains in a compact subset of the state space on $I_{x_0,u} = [0,t_{x_0,u})$. We again conclude by \cite[Prop~C.3.6]{Sontag2013} that $I_{x_0,u} = [0,\infty)$.
	\end{proof}
    
    Observe that, the choice of the control space $\U([0,\infty),\mbbR^m)$ results in a different sufficient condition in Theorem~\ref{thm:IVP_main}.
	In practical applications, the inequalities (A2.1) and (A2.2) may hold only after applying a constant shift to the exergy function. The following remark addresses this situation. 
	\begin{remark}\label{rem:IVP_main}
		The statement of Theorem~\ref{thm:IVP_main} remains true upon a constant shift of the exergy function, i.e., replacing $E(x)$ by $\tilde E(x) := E(x) + D$ for any $D\in \mathbb{R}$. Whereas (A1) holds for $\tilde E$ iff (A1) holds for $E$, this change in particular renders (A2) suitable also for exergy functions bounded from below with a negative minimum. 
	\end{remark}
	
	\subsection*{Discussion on the choice of $E$ as a Lyapunov-like function in Theorem~\ref{thm:IVP_main}}
	
	The essence of the proof of Theorem~\ref{thm:IVP_main} is to show that the unique solution, $x$, to $\IVP$ is constrained to a compact set over its maximal interval of existence, $I_{x_0,u}$, to then conclude via the well-established tool from ordinary differential equations that the solution exists globally. 
	This methodology is closely related to Lyapunov stability theory, see for example \cite[Ch.~4]{Khalil2002}, where one seeks for a function of the state that is \emph{proper} (a continuous function $V:\mbbR^n \rightarrow \mbbR$ is proper if the pre-image of every compact set is compact, that is $V^{-1}(K)$ is compact for every compact $K\subseteq\mbbR^n$) and that remains bounded along any solution of $\IVP$, so that one may conclude that the solution remains in a compact set. 
	Here, leveraging the structure of RI-PHSs, first natural candidates for such a proper function are the Hamiltonian corresponding to the total energy, and the total entropy, that is, $H$ and $S$, respectively.  
	
	However, for irreversible processes, the Hamiltonian $H$ by itself is not a good choice for such a function as for thermodynamic systems, such as the heat exchanger and the gas-piston example introduced in Section~\ref{sec:pH_systems}, it often contains exponentials, and as a consequence, it is not proper. 
	Moreover, as the entropy variables of the subsystems are often included as state variables (see again Section~\ref{sec:pH_systems}) the total entropy $S$ is, as a sum of these states, a linear function. Thus, while using the entropy balance \eqref{eq:entropy_balance} one could deduce a decrease of the negentropy $-S$, this decrease as such does not imply boundedness of solutions, as the total entropy, being a linear function of the states, is not bounded from below.
	
	A further alternative might be the \emph{energy based availability function}, $A:\mbbR^n\rightarrow \mbbR_{\geq 0}$, see \cite{ramirez2016passivity},
	\begin{equation} \label{eq:AvailabilityFunction}
	A(x,x^{\mathrm{eq}}) := H(x) - H(x^{\mathrm{eq}}) - H_x(x^{\mathrm{eq}})^\top (x - x^{\mathrm{eq}}),
	\end{equation}
	where $x^{\mathrm{eq}}\in\mbbR^n$ is a controlled equilibrium point (i.e., a point for which there exists a $u^{\mathrm{eq}}\in\mbbR^m$ such that $f(x^{\mathrm{eq}}, u^{\mathrm{eq}}) = 0$) for which $A(x^{\mathrm{eq}},x^{\mathrm{eq}}) = 0$. 
	While this function is very useful in feedback design, that is, designing \emph{particular inputs} that render $x^{\mathrm{eq}}$ locally asymptotically stable, it is not applicable for addressing the general question of existence of solutions to $\IVP$s for an \emph{a priori given} input.
	
	
	Here, the exergy $E(x) = H(x) - T_0S(x)$ turns out to be a good choice 
	as $H$, involving exponentials of the state, tends to infinity in any direction with positive components whereas the $-S$ term blows up in any direction with purely negative components. More precisely, we will show that for the heat exchanger network, $E$ is \emph{radially unbounded}, that is, $\| x \| \rightarrow \infty$ implies $E(x) \rightarrow \infty$ and thus proper. 
	For the gas-piston system, we will provide a relaxation of the radial unboundedness condition. Therein, we require $E$ to only blow up along \emph{possible} unbounded directions, i.e., (A1). 
	
	\subsection{A corollary for dissipative port-Hamiltonian systems}\label{subsec:DPHS_corollary}
	
	We now show that the result in Theorem~\ref{thm:IVP_main} also applies to \emph{dissipative} port-Hamiltonian (DPH) systems. 
	These are systems that take the form
    \begin{subequations}\label{eq:dph_ges}
	\begin{align}
	\dot x & = \Bigl(J(x)-R(x)\Bigr)H_x(x) + g(x, H_x)u \label{eq:dph}\\
		y_H  &= g(x,H_x)^\top H_x,
	\end{align}
    \end{subequations}
    where $J:\mbbR^n\rightarrow \mbbR^{n \times n}$ is pointwise skew-symmetric and locally Lipschitz, $R:\mbbR^n\rightarrow \mbbR^{n\times n}$ is pointwise positive definite and symmetric, and $g:\mbbR^n \times \mbbR^n \rightarrow \mbbR^{n\times m}$ is locally Lipschitz\footnote{We note that in the usual definition of dissipative port-Hamiltonian systems, cf.\ \cite{van2014port} and the references therein, the input matrix only depends on the intensive and not on the extensive variables, i.e., $g(x,H_x) \equiv g(x)$. Here, we also allow an explicit dependence on $H_x$ in view of Definition~\ref{def:RIPHS}.}.  
	As before, 
	$H:\mbbR^n\rightarrow \mbbR$ is the differentiable Hamiltonian function, which satisfies the energy balance equation
	\begin{align}\label{eq:dissipation_ineq}
	\frac{\dee H}{\dee t}(x(t)) = H_x^\top R(x) H_x + u(t)^\top y_H(t) \leq u(t)^\top y_H(t). 
	\end{align}
	\begin{corollary}\label{cor:dissipative}
		Assume the problem $\IVP$ with 
        a DPH system, i.e., $f(x,u)$ is as in \eqref{eq:dph}. 		If (A1) and (A2) hold with $E$ replaced by $H$, then $I_{x_0,u} = [0,\infty)$. 
	\end{corollary}
	\begin{proof}
		The proof follows directly by the dissipation inequality \eqref{eq:dissipation_ineq} along the same lines of the proof of Theorem~\ref{thm:IVP_main}. 
	\end{proof}
	
	As an application of this result, we show that its assumptions are satisfied e.g.~for quadratic and positive definite Hamiltonians with bounded input matrix. 
	\begin{proposition}\label{prop:DPHS_quad}
		Assume a dissipative Hamiltonian system \eqref{eq:dph_ges} with
		\begin{itemize}
			\item[1)] quadratic Hamiltonian, i.e., $H(x) = \frac{1}{2}\| x \|_Q^2$, $Q\in\mbbR^{n\times n}$, $Q \succ 0$,
			\item[2)] bounded input matrix, i.e., there exists a $K\geq 0$ such that 
			$\sup_{x\in\mbbR^n}\|g(x)\|_M = K,$
			where $\| \boldsymbol\cdot \|_M$ indicates an arbitrary induced matrix norm,
                \item[3)] $u \in \Lcal^1([0,\infty), \mbbR^m)$ 		\end{itemize}
		Then (A1) and (A2.1) are satisfied.
	\end{proposition}
	\begin{proof}
		Clearly, positive definite quadratic Hamiltonians as above are radially unbounded, thus (A1) holds with $E$ replaced by $H$. 
		Moreover, $
		\| H_x (x)\|^2 = \| Qx \|^2 = x^\top Q^\top Q x = \| x \|_{\hat Q}^2,$
		where $\hat Q := Q^\top Q \succ 0$. 
		Using,
		\[
		\tfrac{1}{2}\underline \sigma(Q) \| x\|^2 \leq H(x) \leq \tfrac{1}{2}\overline \sigma(Q)\| x\|^2\quad \text{and}\quad \underline \sigma(\hat Q) \| x\|^2 \leq 	\| H_x(x) \|^2 \leq \overline \sigma(\hat Q)\| x\|^2, 
		\]
		we conclude that
$\| g(x)^\top H_x(x) \|^2 \leq K^2 \|H_x(x) \|^2 \leq 2 K^2 \tfrac{\overline \sigma(\hat Q)}{\underline \sigma(Q)} H(x) $
		for all $x \in \mbbR^n$.
		Thus, with $u\in\Lcal^1([0,\infty),\mbbR^m)$, (A2.1) holds with $E$ replaced with $H$, where $p = 2$ and $C = 2K^2 \tfrac{\overline \sigma(\hat Q)}{\underline \sigma(Q)}$.  
	\end{proof}
    We briefly comment on this result.
	\begin{remark}
            \hspace*{.8em} (i) In Proposition~\ref{prop:DPHS_quad} the conditions 1), 2) and $u\in\Lcal^{\infty}([0,\infty), \mbbR^m)$ do \emph{not} imply that (A2.2) holds.  
            By similar arguments as those in the above proof it can be shown that
            $
                \frac{\|g(x,H_x)^TH_x(x)\|}{H(x)}\leq \frac{K\|Q\|_M\sqrt{\tfrac{2}{\underline{\sigma}(Q)}}}{\sqrt{H(x)}},
            $
             and so the left-hand side only remains finite for $x$ away from the origin.\\[.5mm]
         \hspace*{.8em} (ii) For linear PHS, i.e., where $J(x) \equiv J$, $R(x) \equiv R$, $g(x) \equiv B$ and where $H(x)$ is quadratic in \eqref{eq:dph_ges}, we should note that the result of \cref{thm:IVP_main} may be of very limited use and in particular is not applicable using $E=H$. To see this, we compute $g(x,H_x(x))^\top H_x = B^\top Q x$ which, as a linear function, can not be bounded from above by the quadratic function $H$ around the origin. Hence, neither (A2.1) nor (A2.2) are satisfied.
            In this setting, however, we stress that global existence directly follows from linear theory, that is, the variation of constants formula. 
        \end{remark}	

\section{{Reversible-irreversible pH Systems}: Existence of optimal solutions}
    \label{sec:OCP_existence}
    
	\noindent Following \cite{philipp2023optimal}, we now consider the following nonlinear state- and control-constrained optimal control problem on a time horizon $T >0$
	\begin{numcases}{\OCP_T:\quad}
	\min\limits_{u\in \mathcal{U}_T, x\in\AC_T} & $\quad \mathcal{I}(u,x) := \int_0^T \ell(x(s),u(s))\  \dee s$, \nonumber \\
	\mathrm{subject\  to:}	& $\quad \dot{x}(t) = f(x(t),u(t)),\quad \mathrm{a.e.\quad} t\in[0,T],$ \label{OCP_eq_1}\\
	& $\quad x(0) = x_0,$\\
	& $\quad x(T) \in \mathbb{X}_T.$ \label{OCP_eq_2}
	\end{numcases}
	Here, $f$ is the right-hand side of the reversible-irreversible port-Hamiltonian system~\eqref{eq:RIPH}, and with weights $\alpha_i \geq 0$ for $i\in[1:3]$, the running cost $\ell:\mbbR^n \times \mbbR^m \rightarrow \mbbR$ is given by
	\begin{align}
	\ell(x,u) & = [\alpha_1 y_H - \alpha_2T_0 y_S]^\top u + \alpha_3 \| C x- y^{\mathrm{ref}} \|^2. \label{eq:special_cost}
	\end{align}
	Here, $T_0 >0$ is a reference value, $y_H$ and $y_S$ are the energy- resp.\ entropy-conjugate outputs defined in \eqref{eq:def_ys}, $C\in\mbbR^{p \times n}$ is an output matrix and $y^{\mathrm{ref}}\in\mbbR^p$ is a vector of output reference values. 
	Thus, $\alpha_1$ and $\alpha_2$ determine the weighting in the cost functional between supplied energy encoded by $y_H^\top u$ (cf.~\eqref{eq:energy_balance}) and entropy growth given by $y_S^\top u$ (cf.~\eqref{eq:entropy_balance}), whereas $\alpha_3$ determines the importance of getting certain outputs to be close to desired reference values. 
	
	Further, $x_0\in\mbbR^n$ is an initial state and $\mathbb{X}_T\subseteq\mbbR^n$ denotes the target set. 
	We let,
	\[
	\mathcal{U}_T := \{u \in\U([0,T],\mbbR^m): u(t) \in \mathbb{U}, \mathrm{\quad a.e. \quad} t\in[0,T] \},
	\]
	where $\mathbb{U}\subseteq \mbbR^m$ is the control constraint set and, as in Section~\ref{sec:IVP}, $\U([0,T],\mbbR^m)$ is either $\Lcal^1([0,T],\mbbR^m)$ or $\Lcal^{\infty}([0,T],\mbbR^m)$. In view of the optimal control problem $\OCPT$, we define the following notions, where we abbreviate $\AC_T = \AC([0,T],\mbbR^n)$.
	
	\begin{definition}
		A pair $(u,x)\in\U_T\times\AC_T$, is \emph{admissible} for $\OCPT$ provided it satisfies \eqref{OCP_eq_1}-\eqref{OCP_eq_2}. 
		We indicate all admissible pairs by,
		$\Omega := \{(u,x)\in\U_T\times\AC_T : \eqref{OCP_eq_1}-\eqref{OCP_eq_2} \,\,\mathrm{hold}  \}. $
		A \emph{solution} to $\OCPT$ is a pair $(u^\star,x^\star)\in\Omega$ that globally minimises the cost functional $\mathcal{I}(u,x)$ over all admissible pairs, i.e., $\mathcal{I}(x^\star, u^\star) \leq \mathcal{I}(u, x)$ for all $(u,x) \in\Omega$. 
	\end{definition} 
	In view of the proof of Theorem~\ref{thm:IVP_main} in Section~\ref{sec:IVP}, we may conclude under Assumptions (A1) and (A2) from Theorem~\ref{thm:IVP_main}, that the state of every admissible pair $(u,x)\in\Omega$ satisfies 
	\begin{equation}
	\|x(t)\| \leq \hat K < \infty,\quad \forall t\in[0,T]. \label{def:K_hat}
	\end{equation}
	for a constant $\hat K\in \mathbb{R}$.
	We will use this fact in the next theorem, which is the second main contribution of this work.
	\begin{theorem}\label{thm:OCP_existence}
		Consider $\OCPT$ and assume that (A1) and (A2) of Theorem~\ref{thm:IVP_main} hold, along with the following additional assumptions:
		\begin{enumerate}
			\item[(A3)] The set $\mathbb{U}$ is compact and $\mathbb{X}_T$ is closed. 
			\item[(A4)] There exists an admissible pair for $\OCPT$, that is, $\Omega \neq \emptyset$. 
		\end{enumerate}
		Then there exists a solution to $\OCPT$. 
	\end{theorem}
	\begin{proof}
		We show that (A1)-(A4) imply the sufficient conditions of Filippov's existence Theorem provided in Theorem \ref{thm:filippov} of the appendix.
		
		Under (A1)-(A2) we know, from Theorem~\ref{thm:IVP_main}, that for any admissible pair the state is constrained to a compact set for all $t\in[0,T]$, as in \eqref{def:K_hat}. 
		Under (A3) the set $\{x_0\} \times \mathbb X _T$ is closed and the set $\{x : \|x\| \leq \hat K\} \times \mathbb U$ is compact. 
		Moreover, the dynamics $f$ and running cost $\ell$ are continuous on $\{x : \|x\| \leq \hat K\} \times \mathbb U$. 
		
		Because the dynamics and the running cost are affine in the control 
		we immediately see that the set of points,
		$
		\{(z^0,z^\top)^\top : z = f(x,u), z^0 = \ell(x,u),\quad u\in\mathbb{U}\},
		$
		lies on a convex plane in $\mbbR^{n + 1}$ for every $x\in\mbbR^n$. Thus, together with (A4), we can deduce the existence of an optimal pair $(u^\star, x^\star)\in\Omega$ from Theorem  \ref{thm:filippov}. 
	\end{proof}
	\begin{remark}\label{remark:free_target}
		If the final state in $\OCPT$ is free, that is, if $\mathbb{X}_T = \mbbR^n$, then via Theorem~\ref{thm:IVP_main}, (A1)-(A2) from Theorem~\ref{thm:IVP_main} imply (A4) from Theorem~\ref{thm:OCP_existence}. 
	\end{remark}

\noindent We now apply the theory for global existence of solutions to initial value problems (Section~\ref{sec:IVP}) and existence of solutions to optimal control problems (Section~\ref{sec:OCP_existence}) to a heat exchanger network and the gas-piston system, i.e., the two most-prominent representatives of the system class. 
	To this end, we verify the conditions of Theorem~\ref{thm:IVP_main} and, thus, demonstrate its applicability. 

    \section{Existence of (optimal) solutions for a network of heat exchangers} \label{sec:app:heat_ex}
\begin{figure}[h]
		\centering
		\includegraphics[width=.5\columnwidth]{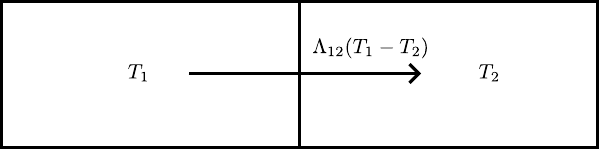}
		\caption{A simple heat exchanger.}\label{fig:two_compartments}
	\end{figure}

    In this section, we apply the existence results of Theorem~\ref{thm:IVP_main} and Theorem~\ref{thm:OCP_existence} to a network of heat exchangers.
	Following \cite[Ch.~2]{ramirez2012control}, 
	we consider two compartments interacting through a heat conducting, non-deformable and impermeable wall as in Figure~\ref{fig:two_compartments}. 
	
	\noindent
	Let $T_1, T_2 \geq 0$ denote the temperature in the first and second compartment, respectively. 
	By Fourier's law we have,
	\begin{equation}
	q  =\Lambda_{12}(T_1 - T_2), \label{eq:Fouriers_law}
	\end{equation}
	where $q$ is the heat flux through the wall
	and $\Lambda_{12}> 0$ is the material's thermal conductivity.  
	The temperature as a function of the entropy is given by
	\begin{equation}
	T_i(S_i) = T_{\mathrm{ref}} e^{\frac{(S_i - S_{\mathrm{ref}})}{c_i}}, \label{eq:temp_v_S}
	\end{equation}
	$i \in [1:2]$, where $T_{\mathrm{ref}}$ and $S_{\mathrm{ref}}$ are a reference temperature and entropy and $c_1, c_2$ are heat capacities, see \cite[Sec.~2.2]{couenne2006bond}. 
	For simplicity, we take $c_1=c_2 = 1$ and $S_{\mathrm{ref}} = 0$. 
	The state of the system is given by the vector of entropies of each compartment, that is $x := (x_1,x_2)^\top = (S_1, S_2)^\top$. 
	The total energy $H$ (resp.\ total entropy $S$) of the system is the sum of the internal energy $H_i$ (resp. entropy $S_i$), $i=1,2$ in each compartment. 
	Therefore, by Gibbs' equation, $\dee H_i = T_i \dee S_i$, the 
	total energy and entropy reads
	\[
	H(x) = H_1(x_1) + H_2(x_2) =T_{\mathrm{ref}}  e^{x_1} + T_{\mathrm{ref}} e^{x_2}, \qquad S(x) = S_1 + S_2 = x_1 + x_2.
	\]
	The change in energy in each compartment equals the heat flux, hence
	\begin{equation}
		q = -\frac{\dee H_1}{\dee t}(x(t)) = -T_1(x_1(t))\dot x_1(t) = \frac{\dee H_2}{\dee t}(x(t)) = T_2 (x_2(t)) \dot x_2(t). \label{eq:q_dHdt}
	\end{equation}
	Using \eqref{eq:Fouriers_law}, \eqref{eq:temp_v_S} and \eqref{eq:q_dHdt}, we see that
	\begin{equation*}
	\dot x_1 = -\Lambda_{12}\tfrac{e^{x_1} - e^{x_2}}{e^{x_1}}, \qquad  \dot x_2    = -\Lambda_{12}\tfrac{e^{x_2} - e^{x_1}}{e^{x_2}}.
	\end{equation*}
	Note that the reference temperature $T_{\mathrm{ref}}$ does not influence the dynamics, as it cancels out in the above fractions. 
    
\begin{figure}[t]
		\centering
		\includegraphics[width= \linewidth]{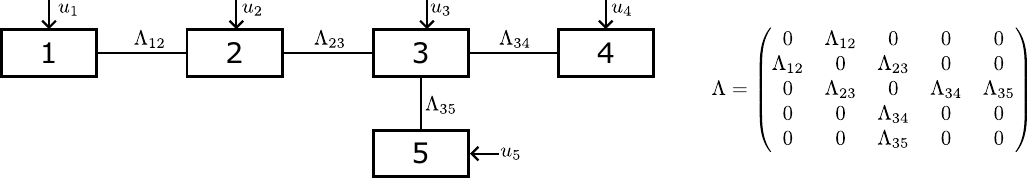}
		\caption{A simple heat exchanger network with thermal conductivity between compartment $i$ and $j$ indicated by $\Lambda_{ij} = \Lambda_{ji}$ and corresponding incidence matrix.}
		\label{fig:heat_exchanger_network}
	\end{figure}

    Consider now a network of $n\in\mathbb N$ coupled compartments, its topology described via an edge-weighted graph, $G(V, \Lambda)$, where $V :=\{V_1,V_2,\dots,V_n\}$ denotes the vertices and $\Lambda\in\mbbR^{n \times n}$ is the incidence matrix, see Figure~\ref{fig:heat_exchanger_network} for an example.
	The $ij$-th entry of the incidence matrix, denoted $\Lambda_{ij} > 0$, is the material conductivity between compartment $i$ and $j$. 
	Thus, $\Lambda$ is symmetric with zeros on its main diagonal. 
	By the same reasoning as in the two-compartment case, assuming each compartment has its own controllable entropy input (denoted by $u_i$) corresponding to local heating/cooling, the differential equation associated with the $i$-th compartment reads
	\begin{align}
	\dot x_i(t) & = -\sum\nolimits_{j = 1}^{n} \Lambda_{i j} \left(\frac{e^{x_i} - e^{x_j}}{e^{x_i}} \right) + u_i, \label{eq:heat_exchanger_nerwork}
	\end{align}
	$i\in[1:n]$, where the state is $x := (x_1,x_2,\dots x_n)^\top = (S_1,S_2,\dots,S_n)^\top\in\mbbR^n$. 
	The total energy and entropy functions of the network read
	\begin{equation}
	H(x) = T_{\mathrm{ref}} \sum\nolimits_{i=1}^n e^{x_i}, \qquad S(x) = \sum\nolimits_{i=1}^n x_i.
	\end{equation}
	The system \eqref{eq:heat_exchanger_nerwork} may be written in the form of \eqref{eq:RIPH} as follows.  
	We set $J_0 = 0$, $B = I_n$ and let
	\[
	\gamma_{(i,j)}(x) = \frac{\Lambda_{ij}}{T_{\mathrm{ref}}^2}e^{x_i}e^{x_j},
	\]
	for all $i\in[2:n]$ and $j\in[i+1,n]$ and let $J^{(i,j)} \in\mbbR^{n\times n}$ be given by
	\begin{equation*}
	J^{(i,j)}_{kl}=\begin{cases}
	-1 & \text{if } (k,l)=(i,j), \\
	1 & \text{if } (k,l)=(j,i),
	\end{cases}
	\end{equation*}
	for all $i\in[2:n]$ and $j\in[i+1,n]$ encoding the $N := \frac{n(n-1)}{2}$ interactions between the compartments. 
	The resulting model then may be formulated as a RI-PHS \eqref{eq:RIPH} generated by the Hamiltonian $H$ being the total energy as follows
\begin{align}\label{eq:phs_heatnetwork}
	\dot x = \biggl(\sum\nolimits_{\substack{i\in[2:n] \\ j\in[i+1:n]}}\gamma_{(i,j)}(x)\{S,H\}_{J^{(i,j)}}(x)J^{(i,j)}\biggr)H_x(x) + u,
	\end{align}
	where $u:=(u_1,u_2,\dots,u_n)^\top$. 
	Note that the brackets define the thermodynamical force giving rise to
	the irreversible phenomenon, here the heat flux resulting from temperature differences
	$
	\{S,H\}_{J^{(i,j)}}(x) = T_{\mathrm{ref}}e^{x_i} - T_{\mathrm{ref}}e^{x_j} = T_i - T_j.
	$
	\begin{theorem}[Existence of solutions]\label{thm_exheat}
	For any $x_0\in \mathbb{R}^n$ and $u\in \mathcal{L}^\infty([0,\infty),\mathbb{R}^m)$, the RI-PHS modeling the network of heat exchangers \eqref{eq:phs_heatnetwork} gives rise to a unique global solution, i.e., $I_{x_0,u} = [0,\infty)$.
	\end{theorem}
    \begin{proof}
    We will verify the assumptions of Theorem~\ref{thm:IVP_main} to show global existence. Since the exergy
	\begin{align*}
	E(x) = H(x) - T_0S(x) = \sum\nolimits_{i=1}^n T_{\mathrm{ref}}e^{x_i} - T_0x_i,
	\end{align*}
 is radially unbounded due to the positivity of $T_0$, that is, if $\|x\| \rightarrow \infty$ then $E(x) \rightarrow \infty$. Consequently, Assumption~(A1) is satisfied. 
	We now prove that (A2.2) holds in the sense of Remark \ref{rem:IVP_main}, i.e., the claimed inequality holds for the shifted exergy function. In other words, we show that there exist constants $C,D > 0$ such that
	\[
	\left| \left| g(x(t),H_x(x(t)))^\top E_x(x(t)) \right|\right| = \left| \left| E_x(x(t)) \right|\right|  \leq C(E(x(t)) + D) \coloneqq C\tilde E(x(t))
	\]
	for all $t \in I_{x_0,u}$.
	To this end, we compute the derivatives of the exergy function
	\[
	E_x(x) = 
	\left(\begin{smallmatrix}
	T_{\mathrm{ref}}e^{x_1} - T_0 \\
	T_{\mathrm{ref}}e^{x_2} - T_0 \\
	\vdots \\
	T_{\mathrm{ref}}e^{x_n} - T_0
	\end{smallmatrix}\right),
	\quad
	E_{xx}(x) = T_{\mathrm{ref}}
	\left(\begin{smallmatrix}
	e^{x_1} & 0 & 0& \dots & 0\\
	0 & e^{x_2} & 0 & \dots & 0 \\
	0 & 0 & e^{x_3} & \dots  &0 \\
	\vdots & \vdots & \vdots & \ddots  & 0 \\
	0 & 0 & 0 & \dots & e^{x_n}
	\end{smallmatrix}\right). 
	\]
	Thus, $E$ is strictly convex by the fact that $E_{xx}$ is positive definite for all $x\in\mbbR^n$. Let $\underline x \coloneqq \ln(T_0/T_{\mathrm{ref}} ) \mathbf{1}_n \in\mbbR^n$, where $\mathbf{1}_n$ denotes the constant vector of ones in $\mbbR^n$. 
	As $E_x(\underline x) = 0$ we deduce that $E$ has a global minimum at $\underline x$. 
	Therefore, 
	$$E(x) \geq E(\underline x) = nT_0\left(1- \ln\left(\frac{T_0}{T_{\mathrm{ref}}}\right)\right)$$ for all $x\in\mbbR^n$. 
	Recall that the state of the system consists of the individual entropy variables. The total entropy is given by $S(x)=\sum\nolimits_{i=1}^n x_i$ and thus its gradient is the (constant) vector of ones. Hence $\|S_x(x)\|_1=n$ and 
	\begin{equation}
	T_0\| S_x(x)\|_1 =T_0n \leq E(x) + nT_0\ln\left(\frac{T_0}{T_{\mathrm{ref}}}\right) \label{ineq_S_bound_heat}
	\end{equation}
	for all $x\in\mbbR^n$. 
	We now show that the quotient
	\begin{align}\label{eq:limit}
	\frac{\|H_x(x)\|_1}{E(x)} = \frac{\sum\nolimits_{i=1}^n T_{\mathrm{ref}}e^{x_i}}{\sum\nolimits_{i=1}^n T_{\mathrm{ref}}e^{x_i} - T_0x_i}
	\end{align}
	is bounded over $\mbbR^n$. 
	To that end, let $c_1,c_2,\dots,c_n \in \mbbR$ be constants that define a ray $r:\mbbR\rightarrow\mbbR^n$ parametrised by $s\geq 0$ passing through the origin, i.e.,
	\[
	r(s) = (c_1 s, c_2 s, \dots, c_n s), \quad s\geq 0.
	\]
	Consider the limit of \eqref{eq:limit} along this ray.
	If $c_i \leq 0$ for all $i\in[1:n]$, and $c_i < 0$ for at least one $i\in[1:n]$, then,
	\[
	\lim_{s \rightarrow \infty} \frac{\sum\nolimits_{i=1}^nT_{\mathrm{ref}}e^{c_i s}}{\sum\nolimits_{i=1}^n T_{\mathrm{ref}}e^{c_i s} - T_0c_i s} = 0.
	\]
	Note that we ignore the case where $c_1 = c_2 = \dots = c_n = 0$ because then $r$ is just the origin. 
	If there exists an $i\in[1:n]$ for which $c_i > 0$, the limit takes the indeterminate form infinity over infinity.
	In this case, by L'H\^opital's rule
	\[
	\lim_{s \rightarrow \infty} \frac{\sum\nolimits_{i=1}^n T_{\mathrm{ref}}e^{c_i s}}{\sum\nolimits_{i=1}^n T_{\mathrm{ref}}e^{c_i s} - T_0c_i s} 
	= 
	\lim_{s \rightarrow \infty} \frac{\sum\nolimits_{i=1}^n T_{\mathrm{ref}}c_i e^{c_i s}}{\sum\nolimits_{i=1}^n T_{\mathrm{ref}}c_i e^{c_i s} - T_0c_i}
	= 1.
	\]
	We deduce that the function $\frac{\|H_x(x)\|_1}{E(x)}$ is bounded along any unbounded direction. 
	By continuity, we conclude that there exists a constant $K < \infty$ such that
$\|H_x(x)\|_1 \leq K E(x)$
	for all $x\in\mbbR^n$.
	Recall that $g(x, H_x(x)) = I_n$ for all $x \in \mbbR^n$. Hence,  by \eqref{ineq_S_bound_heat},
	\begin{align*}
	\|E_x(x)\|_1 &\leq \|H_x(x)\|_1 + T_0\|S_x(x)\|_1 \leq (K + 1)\bigg(E(x) + \tfrac{nT_0\ln\left(\tfrac{T_0}{T_{\mathrm{ref}}}\right)}{K + 1}\bigg)
	\end{align*}
	for all $x\in\mbbR^n$. 
 Finally, with $u\in\Lcal^{\infty}([0,\infty),\mbbR^m)$, (A2.2) holds for an arbitrary heat exchanger network with $C =K + 1$ and $D=\tfrac{nT_0\ln\left(T_0/T_{\mathrm{ref}}\right)}{K + 1}$, cf.\ Remark \ref{rem:IVP_main}. This concludes the proof. 
 \end{proof}

Having analyzed the initial value problem, we now may straightforwardly invoke Theorem~\ref{thm:OCP_existence} to deduce existence of solutions to the corresponding entropy-optimal control problem $(\OCP_T)$ with suitable control and terminal constraint sets. 
\begin{corollary}
	For any $T>0$, compact control set $\mathbb{U}\subset \mathbb{R}^m$ and closed terminal set $\mathbb{X}_T$ that is reachable with a measurable control $u:[0,T]\to \mathbb{U}$, the problem $\OCPT$ subject to the network of heat exchangers~\eqref{eq:heat_exchanger_nerwork} possesses a solution. 
\end{corollary}
		
	\section{Existence of (optimal) solutions for a gas-piston System}\label{sec:app:gaspiston}
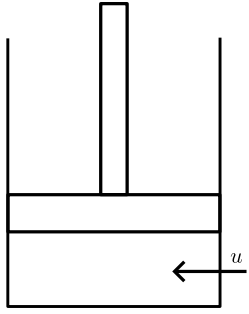
\begin{figure}[h]
		\centering
		\scalebox{0.7}{
        \begin{tikzpicture}[]

  \draw (0,0) -- (0,3.75);
  \draw (4,0) -- (4,3.75);
  \draw (0,0) -- (4,0);

  \draw[fill=gray, fill opacity=0.2, draw=gray, dashed] (0.05,1.5) -- (1.8,1.5) -- (1.8,3.5) -- (2.2,3.5) -- (2.2,1.5) -- (3.95,1.5) -- (3.95,1.3) -- (0.05,1.3) -- cycle;


  \node at (5.0, 0.7)  {$u$};

  \draw[->, thick] (4.8,0.7) -- (3.2,0.7) ;
\end{tikzpicture}}
	\caption{A gas-piston system with control as entropy flow into the cylinder.}\label{fig:gas_piston}
\end{figure}

	Second, we provide an application of our main results to a gas-piston system. Consider a cylinder filled with an ideal gas, closed by a piston subject to gravity, as in Figure~\ref{fig:gas_piston}, see \cite{philipp2023optimal, eberard2007extension}.
We let $S\in\mbbR$ denote the system's entropy, $V\in\mbbR$ the volume of the gas in the cylinder and $p\in\mbbR$ the momentum of the piston. 
	The system's internal energy reads 
	\[
	U(S,V) = \frac{3}{2}NRT_{\mathrm{ref}}\cdot e^{\beta(S,V)}, \qquad \beta(S,V) = \frac{S - Ns_{\mathrm{ref}}}{\frac 32 RN} + \frac{2}{3}\left( \ln(NRT_{\mathrm{ref}}) - \ln(VP_{\mathrm{ref}}) \right),
	\]
	where $N > 0$ is the number of moles of gas in the cylinder, $R>0$ is the ideal gas constant and $s_{\mathrm{ref}}$, $T_{\mathrm{ref}}$, $P_{\mathrm{ref}} > 0$ are positive reference values. 
	The Hamiltonian is the sum of the system's internal, kinetic and potential energies and is given by
	\begin{equation}
	H(S, V, p) = U(S,V) + \frac{1}{2m} p^2 + \frac{mg}{A}V, \label{eq:H_gas_piston}
	\end{equation}
	where $m > 0$ is the piston's mass, $g > 0$ is the gravitational constant and $A >0$ is the area of the piston head. 
	The temperature and the pressure in the cylinder and the piston's velocity, expressed as a function of the entropy and volume, read
	\[
	T(S,V) = \frac{\partial H}{\partial S} = T_{\mathrm{ref}} e^{\beta(S,V)},\ 	P(S,V)   = -\frac{\partial U}{\partial V} =  \frac{NRT_{\mathrm{ref}} e^{\beta(S,V)}}{V},\ 	v(p) = \frac{\partial H}{\partial p} = \frac{1}{m} p
	\]
	respectively.
%
%
	We denote the system's state by $x := (x_1,x_2,x_3)^\top = (S, V, p)$ and we consider the control to be the entropy supplied to the cylinder. 
	Therefore, as explained in detail in \cite{philipp2023optimal} the equations of motion read
	\begin{align}\label{eq:gas_piston_S}
	\frac{\mathrm{d}}{\mathrm{d}t} \begin{pmatrix}
	x_1 \\ x_2 \\ x_3
	\end{pmatrix} &= \begin{pmatrix}
	\frac{\kappa v(x_3)^2}{T(x_1,x_2)} \\
	Av(x_3) \\
	-\kappa v(x_3) + AP(x_1,x_2)
	\end{pmatrix} + \begin{pmatrix}
	u \\ 0 \\ 0
	\end{pmatrix},
	\end{align}
	where $u\in\mbbR$ is the supplied entropy and $\kappa > 0$ is a positive damping constant. 
	Hence, \eqref{eq:gas_piston_S} is a RI-PHS by Definition \ref{def:RIPHS} where
	\[
	g =\left(\begin{smallmatrix}
	1 \\ 0 \\ 0
	\end{smallmatrix}\right), \quad
	J_0 = 
	\left(\begin{smallmatrix}
	0 & 0 & 0 \\
	0 & 0 & A \\
	0 & -A & 0
	\end{smallmatrix}\right),
	\quad 
	J_1 = 
	\left(\begin{smallmatrix}
	0 & 0 & 1 \\
	0 & 0 & 0 \\
	-1 & 0 & 0
	\end{smallmatrix}\right)\quad \text{and} \quad \gamma(x) =  \tfrac{\kappa}{T(x_1,x_2)}.
	\]
	Note that $J_0 S_x = J_0\begin{pmatrix}
	1 & 0 & 0
	\end{pmatrix}^\top = 0$ and so $S$ is a Casimir function of $J_0$.

    \begin{theorem}[Global existence of solutions]
    For any $x_0\in \mathbb{R}^n$ and $u\in \mathcal{L}^\infty([0,\infty),\mathbb{R}^m)$, the RI-PHS modeling the gas-piston system~\eqref{eq:gas_piston_S} gives rise to a unique global solution, i.e., $I_{x_0,u} = [0,\infty)$.
    \end{theorem}
    \begin{proof}
    Again, we will verify the assumptions of Theorem~\ref{thm:IVP_main}. For simplicity, set
	\[
	e^{\beta(x_1,x_2)} = K_1 \tfrac{e^{K_2 x_1}}{x_2^{2/3}},
	\]
	where
	$K_1 := \frac{3}{2}\tfrac{(NRT_{\mathrm{ref}})^{5/3}}{P_{\mathrm{ref}}^{2/3}} \exp\left({-\tfrac{2}{3}\frac{s_{\mathrm{ref}}}{R}}\right), K_2 := \tfrac{2}{3RN}$ are positive constants.
	We first show that the system's state is constrained to the domain $D:=\mbbR\times\mbbR_{>0}\times\mbbR$. 
	To that end, consider an arbitrary initial state $x_0 = x(0) \in D$, with $x_1(0)\in \mbbR$, $x_2(0) > 0$ and $x_3(0) < 0$. 
	Referring to the equations of motion, \eqref{eq:gas_piston_S}, we see that $\dot x_2(0) < 0$ and $\dot x_3(0) > 0$. 
	Thus, for any $u\in\Lcal^{\infty}([0,\infty),\mbbR)$ there exists a nondegenerate interval, $[0,\bar t]$, over which $x_2$ is strictly decreasing and $x_3$ is strictly increasing.  
	Moreover, by continuity of the state trajectory, the constant $K_0 := \min_{t\in[0,\bar t]}e^{K_2 x_1(t)}>0$ is finite. 
	Therefore, referring now to the third equation of \eqref{eq:gas_piston_S} over the interval $[0,\bar t]$, we see that, for any $t\in[0,\bar t],$
	\begin{align*}
	\dot x_3(t) & = -\frac{\kappa}{m}x_3(t) + \frac{ANRT_{\mathrm{ref}} K_1 e^{K_2 x_1(t)}}{(x_2(t))^{5/3}} \geq -\frac{\kappa}{m}x_3(t) + \frac{\hat K_0 }{x_2(0)^{5/3}},
	\end{align*}
	where we let $
	\hat K_0 := ANRT_{\mathrm{ref}} K_1 K_0 > 0.$
	By Gr\"onwall's Lemma and the variation of constants formula we get,
	\[
	x_3(t) \geq e^{-\frac{\kappa}{m}t}x_3(0) + \frac{\hat K_0}{x_2(0)^{5/3}}(1 - e^{-\frac{\kappa}{m}t}),\quad t\in[0,\bar t],
	\]
	and integrating the second equation in \eqref{eq:gas_piston_S} we get, for any $t\in[0,\bar t]$,
	\begin{align*}
	x_2(t)  = x_2(0) \!+\! \frac{A}{m}\int_0^t \!x_3(s)\, \dee s &\geq x_2(0) \!+\! \frac{A}{m}\int_0^t \!\bigg( e^{-\frac{\kappa}{m}s}x_3(0) + \tfrac{\hat K_0}{x_2(0)^{5/3}}(1 - e^{-\frac{\kappa}{m}s})\,\bigg)\dee s \\
	& = x_2(0) + \hat K_1 t + \hat K_2\left( e^{-\frac{\kappa}{m}t} - 1 \right) := x_2(0) + G(t),
	\end{align*}
	where
	$\hat K_1 := \frac{A}{m}\frac{\hat K_0}{x_2(0)^{5/3}}>0, 
	\hat K_2 := \frac{A}{\kappa}\left(\frac{\hat K_0}{x_2(0)^{5/3}}-x_3(0) \right).$
	Consider the function $G(t)$ on $(-\infty,\infty)$. 
		Requiring $\frac{\dee G}{\dee t} = 0$ and recalling $x_2(0)>0$ and $x_3(0) < 0$, we see that $G(t)$ attains an extremal value at
	\[
	\hat t = -\frac{m}{\kappa}\ln\left(\frac{m \hat K_1}{\kappa \hat K_2}\right) = -\frac{m}{\kappa}\ln\bigg(\big(1 - \frac{x_3(0)x_2(0)^{5/3}}{\hat K_0}\big)^{-1}\bigg) > 0. 
	\]
	Moreover, considering $\frac{\dee^2 G}{\dee t^2} > 0$, this extremal value is a global minimum of $G$ (on its domain $(-\infty,\infty)$ ), provided $\hat K_2 > 0$, which in turn holds if $x_2(0)>0$ is sufficiently small. 
	Finally, substituting $\hat t$ into $G$, we see that,
	\[
	G(\hat t) = \hat K_1\left(\hat t + \frac{m}{\kappa}\right) - \hat K_2 = \frac{A}{m}\frac{\hat K_0}{x_2(0)^{5/3}}\hat t + \frac{Ax_3(0)}{\kappa},
	\]
	which is also strictly positive provided $x_2(0)$ is sufficiently small. 
	Thus, with $x_2(0)>0$ sufficiently small, $G(t) > 0$ for all $t\geq 0$, and so $x_2(t) > 0$ for all $t\in[0,\bar t]$. 
	Now, with $x_1(0)\in\mbbR$, $x_2(0)>0$ but with $x_3(0)\geq 0$, one immediately sees from $\dot x_2(0) \geq 0$ that $x_2(t)>0$ for all $t$ in a nondegenerate interval. 
	Because the dynamics, $f$, is time-invariant, we can conclude that for any $x_0\in D$ and any $u\in\Lcal^{\infty}([0,\infty),\mbbR)$, $x_2(t) > 0$ for all $t\in I_{x_0,u}$, its maximal interval of existence. 
	
	We now show that (A1) holds. 
	To that end, note that we may express the exergy function as follows,
	\begin{align}
	E(x) & = H(x) - T_0S(x) 
	= K_1 \tfrac{e^{K_2 x_1}}{x_2^{2/3}} + K_3x_3^2 + K_4 x_2 - T_0x_1, \label{eq:frak_S_Gas_piston}
	\end{align}
	where $K_1$ and $K_2$ are as above, and $K_3 := \frac{1}{2m}, K_4 := \frac{mg}{A}$ are positive constants.
	Consider $u~\in~\Lcal^{\infty}([0,\infty),\mbbR)$ and an unbounded local solution, that is, one for which,
	$\lim_{t\rightarrow t_{x_0,u}} \| x(t) \| = \infty.$
	As
	\begin{equation}
	\frac{\partial E}{\partial x_1} = \frac{K_1 K_2}{x_2^{2/3}}e^{K_2 x_1} - T_0, \label{eq:G_partial}
	\end{equation}
	we see that for fixed $x_2> 0$ (corresponding to positive volume) and $x_3\in \mathbb{R}$ the function $E$ has a global minimum at
	\begin{equation}
	\hat x_1 = \frac{1}{K_2}\ln\left(\frac{T_0x_2^{2/3}}{K_1K_2}\right). \label{eq:hat_x_1}
	\end{equation}
	Substituting $\hat x_1$ into \eqref{eq:frak_S_Gas_piston} we get
	\begin{equation}
	E(x) \geq \frac{T_0}{K_2} + K_3x_3^2 + K_4x_2 - \frac{T_0}{K_2}\ln\left( \frac{T_0x_2^{2/3}}{K_1 K_2} \right) \label{ineq:E_lower:bound}
	\end{equation}
	for all $x=(x_1,x_2,x_3)\in \mbbR^3$ with $ x_2> 0$. 
	Therefore, if along an unbounded state trajectory $x(t)=(x_1(t),x_2(t),x_3(t))$ we have $x_2(t)\rightarrow \infty$ or $x_3(t) \rightarrow \pm \infty$ then  $\lim_{t\rightarrow \infty} E(x(t))  = \infty$ via \eqref{ineq:E_lower:bound}.
	
	It remains to show that $E(t)\rightarrow \infty$ for trajectories with $x_1(t)\rightarrow \pm \infty$, but for which $x_2(t)$ and $x_3(t)$ remain bounded.
	To this end, consider a sequence of points $(x^i)_{i\in \mathbb{N}} := \big((x_1^i, x_2^i, x_3^i)\big)_{i\in \mathbb{N}}$ along the unbounded solution $x(t)$ with $x_2^i\leq \overline x_2 <\infty$, $x_3^i \leq \overline x_3 < \infty$ but with $x_1^i \rightarrow \infty$. 
	Due to \eqref{eq:G_partial} and \eqref{eq:hat_x_1}, $\frac{\partial E}{\partial x_1}(x) =\frac{\partial E}{\partial x_1}(x_1)  > 0$ if and only if $x_1 > \hat x_1$. 
	Thus, because $x_2^i$ is bounded there is $\bar i\in \mathbb{N}$ such that $x_1^i > \frac{1}{K_2}\ln(\frac{ T_0 \overline x_2^{2/3}}{K_1K_2}) \geq \frac{1}{K_2}\ln(\frac{T_0(x_2^i)^{2/3}}{K_1K_2})$ for all $i \geq \bar i$, i.e.\ for all $i \geq \bar i$, $E$ is strictly increasing and
	\[
	\lim_{i\rightarrow \infty} E(x^i) = \infty.
	\]
	Similarly, consider now a sequence of points $x^i$ along the unbounded solution $x(t)$, with $x_2^i\leq \overline x_2 <\infty$, $x_3^i \leq \overline x_3 < \infty$ but with $x_1^i \rightarrow -\infty$. 
	As the first state corresponds to the entropy variable, we briefly note that the decrease of entropy is possible because of the presence of the control $u$ in \eqref{eq:gas_piston_S}, which can be negative, i.e., corresponding to cooling. 
	Directly from \eqref{eq:frak_S_Gas_piston} we see that 
	$E(x^i) \xrightarrow{i\rightarrow\infty} \infty$ 
	such that (A1) holds. 
	
	We now proceed to show (A2.2). As $E$ is bounded from below, we may set
	\begin{align*}
	\tilde E(x) := E(x) - \min\bigg(\min_{x\in \mathbb{R}^n} E(x),0\bigg) + \delta
	\end{align*}
	with arbitrary but fixed $\delta > 0$ such that $\tilde E(x) \geq E(x)$ for all $x\in \mathbb{R}^n$ and
	$\tilde{E}(x) \geq \delta > 0.$
	In view of Remark~\ref{rem:IVP_main}, it is also sufficient to show (A2.2) for $\tilde E$ instead of $E$, where we note that $E_x = \tilde E_x$.
	Recall that $g = \begin{pmatrix}
	1 & 0 & 0
	\end{pmatrix}^\top$. 
	Therefore, using the definition of $E$ given in \eqref{eq:frak_S_Gas_piston}, we see that
	\begin{align*}
	g(x, H_x(x))^\top E_x(x) & =  \frac{K_1 K_2}{x_2^{2/3}}e^{K_2 x_1} -  T_0 \\
	& \leq K_2E(x) - K_2\left( K_3 x_3^2 + K_4 x_2 -  T_0x_1 \right) \leq  K_2 (E(x) +  T_0x_1)
	\end{align*}
	for all $ x_2> 0, x_3\in\mbbR$.
	We immediately get that 
	\[
	g(x, H_x(x))^\top E_x(x) \leq K_2 E(x)
	\]
	for all $x_1 \leq 0, x_2 > 0, x_3\in\mbbR$. We still remain to consider the case where $x_1 > 0$. 
	To this end we may estimate
	\[
	g(x, H_x(x))^\top E_x(x) \leq K_2 (E(x) +  T_0x_1) \leq K_2 (\tilde E(x) +  T_0x_1) = K_2 \tilde{E}(x)\left(1 +  \frac{ T_0x_1}{\tilde E(x)} \right)
	\]
	for all $ x_1 >0,  x_2 > 0, x_3\in\mbbR$. Let us now consider the case $x_1\to 0$, that is, we may w.l.o.g.\ assume that $x_1$ is bounded. In this case,
	\[
	g(x, H_x(x))^\top E_x(x) \leq \tilde E(x)\left(1 +  \frac{ T_0x_1}{\tilde E(x)} \right) \leq \tilde E(x)\left(1 + \frac{T_0x_1}{\delta}\right) \leq C\tilde E(x),
	\]
	for a constant $C\geq 0$ which yields the claim.
	
	We remain to show that $T_0x_1/\tilde E(x)$ remains bounded also in the case $x_1\rightarrow \infty$. Due to radial unboundedness $\tilde{E}(x)\rightarrow \infty$.
	Accordingly, as with the heat exchanger network, consider the ray originating at the origin
	$r(s) =   \begin{pmatrix}
	c_1s & c_2s & c_3s
	\end{pmatrix} = \begin{pmatrix}
	x_1 & x_2 & x_3
	\end{pmatrix} $
	with $ s\geq 0$.
	Here we let $c_1 > 0$ since we are only interested in positive $x_1$. Further, we let $c_2>0$ as the volume must be positive and $c_3 \in\mbbR$ because the momentum may be positive or negative.  
	Now we consider the limit, in which the denominator is bounded from below by definition
	\begin{align*}
	\lim_{s \rightarrow \infty} \frac{T_0 x_1}{\tilde E(x)} = \lim_{s\rightarrow\infty} \frac{c_1 s}{\frac{K_1}{(c_2 s)^{2/3}}e^{K_2 c_1 s} + K_3(c_3 s)^2 + K_4 c_2 s - T_0c_1 s } = 0.
	\end{align*}
	By continuity we may deduce that there exists a $C_1 < \infty$ such that
	$\frac{T_0x_1}{E(x)} \leq C_1$	for all $ x\in D$.
	Thus, we have shown that
	\[
	\|g(x, H_x(x))^\top E_x(x)\|_1 \leq K_2 \left(1 + C_1 \right) E(x)
	\]
	for all $x_1 \in\mbbR, x_2 >  0, x_3 \in \mbbR$. So, with $u\in\Lcal^{\infty}([0,\infty), \mbbR^m)$, (A2.2) holds for the gas-piston system with $C = K_2 \left(1 + C_1 \right)$. Thus, we may conclude global existence of solutions via Theorem~\ref{thm:IVP_main}.
    \end{proof}

	In view of entropy- and energy-optimal control, one may now straightforwardly deduce the existence of global solutions for the Problem~$\OCPT$.
    
    \begin{corollary}
        For any $T>0$, compact control set $\mathbb{U}\subset \mathbb{R}^m$ and closed terminal set $\mathbb{X}_T$ that is reachable with a measurable control $u:[0,T]\to \mathbb{U}$, the problem $\OCPT$ subject to the gas-piston system~\eqref{eq:gas_piston_S} possesses a solution. 
    \end{corollary}
		\section{Numerical examples}\label{sec:numerics}
	
	\noindent We now present some numerical experiments to point out possible directions of future research. To this end, we inspect the heat exchanger network without terminal constraints. For numerical examples involving the gas piston system and terminal conditions, we refer to \cite{philipp2023optimal}.
    	
	\subsection{Optimal control of heat exchanger network}\label{subsec:numerics_OCP}
	Consider the optimal control of the five-compartment heat exchanger network shown in Figure~\ref{fig:heat_exchanger_network}, introduced in \cite{philipp2023optimal}, 
	\begin{numcases}{(\OCP_T^{\mathrm{heat}}):}
	\min\limits_{u\in\U_T, x\in\AC_T} &$\mathcal{I}(u, x) := \int_0^T \ell(x(s), u(s))\  \dee s$, \nonumber \\
	\mathrm{subject\  to:}	&$\dot x_i(t)  \!=\! -\sum\nolimits_{j=1}^5 \Lambda_{ij}\!\left(\frac{e^{x_i} - e^{x_j}}{e^{x_i}}\right) + u_i(t),\, \mathrm{a.e.\,} t\in[0,T]$,\nonumber\\
	& $\ x(0) = x_0,$ \nonumber
	\end{numcases}
	with $T_{\mathrm{ref}} = c_i = 1$ for $i\in[1:5]$ and $S_{\mathrm{ref}} = 0$.  
	We consider the running cost \eqref{eq:special_cost} with $T_0 = 1$, that is
	\[
	\ell(x,u)  = [\alpha_1 y_H - \alpha_2T_0 y_S]^\top u + \alpha_3 \| C x- y^{\mathrm{ref}} \|^2
	\]
	with $C = 
	\left(\begin{smallmatrix}
	1 & 0 & 0 & 0 & 0 \\
	0 & 0 & 1 & 0 & 0 \\
	0 & 0 & 0 & 1 & 0
	\end{smallmatrix}\right)$ and $
	y_{\mathrm{ref}} = 
	\left(\begin{smallmatrix}
	1 &
	5 &
	10
	\end{smallmatrix}\right)^\top$.
	Moreover, we take the control constraint to be $\mathbb{U} = [-\bar u, \bar u]\times\{0\}\times \{0\}\times \{0\} \times [-\bar u, \bar u]$, $\bar u > 0$, which implies that we can only control compartments 1 and 5. 
	The cost functional in the optimal control problem now models the control task of driving the entropy in the compartments $1$,\ $3$ and $4$ close to the reference values specified in $y_{\mathrm{ref}}$ via $\| C x- y^{\mathrm{ref}} \|^2$, while penalising entropy withdrawn from the system (via $-y_S^\top u$) and energy supplied into the system via $y_H^\top u$. Further, we may reformulate the cost function via \eqref{eq:energy_balance} and \eqref{eq:entropy_balance} and observe that it is indeed composed by positive semi-definite integral terms and terminal costs
	\begin{align}
	\mathcal{I}(u,x) & = \int_0^T \ell(x(s), u(s)) \,\, \dee s \label{eq:integral_cost}\\
	& = \alpha_1\biggl[H(x(T)) - H(x_0)\biggr] + \alpha_2T_0 \biggl[ S(x_0) - S(x(T)) \biggr] \nonumber\\ 
	& +  \int_0^T\biggl(\alpha_2T_0\biggl[\sum\nolimits_{k=1}^N \gamma_k(s) \left(\{S, H\}_{J_k}x(s)\right)^2\biggr]  + \alpha_3 \| C x(s)- y^{\mathrm{ref}} \|^2\biggr)\, \dee s. \nonumber
	\end{align}
	where we abbreviated $\gamma_k(s) := \gamma_k(x(s), H_x(x(s)))$.
	Given $x_0\in\mbbR^5$ and $u\in\U_T$, let $x(t;x_0,u)$ denote the unique solution to the differential equation at time $t\geq 0$ obtained with control $u$, initiating from $x_0$. 
	
\begin{figure}[htb] 
		\centering
		\begin{subfigure}[t]{0.45\linewidth}
			\centering
			\includegraphics[width=\linewidth]{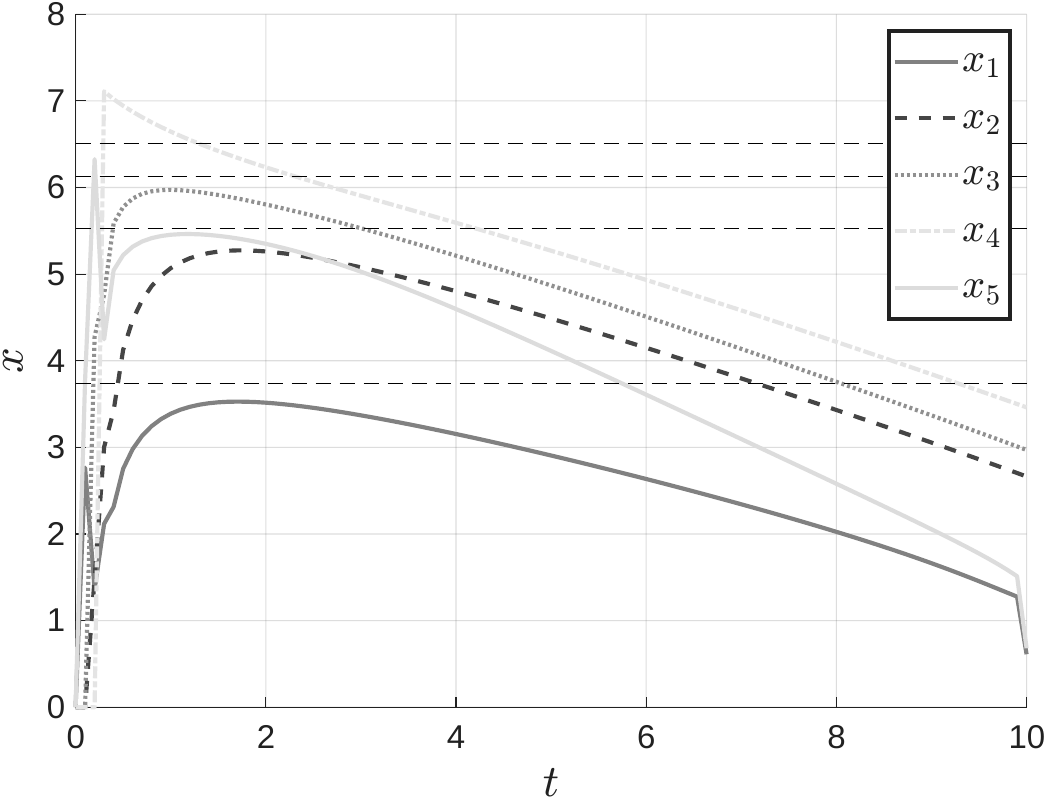}
			\caption{$T = 10$}
			\label{fig:OCP_experiments_a}
		\end{subfigure}\hfil
		\begin{subfigure}[t]{0.45\linewidth}
			\centering
			\includegraphics[width=\linewidth]{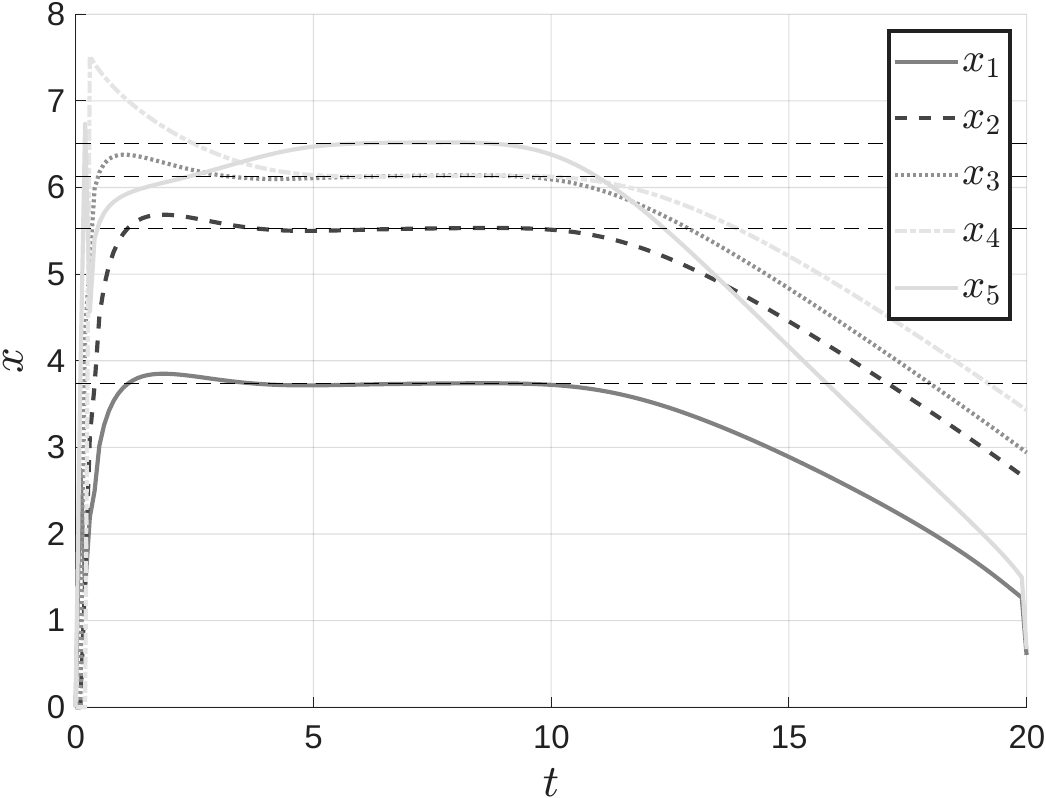}
			\caption[b]{$T = 20$}
		\end{subfigure}
		\begin{subfigure}[t]{0.45\linewidth}
			\centering
			\includegraphics[width=\linewidth]{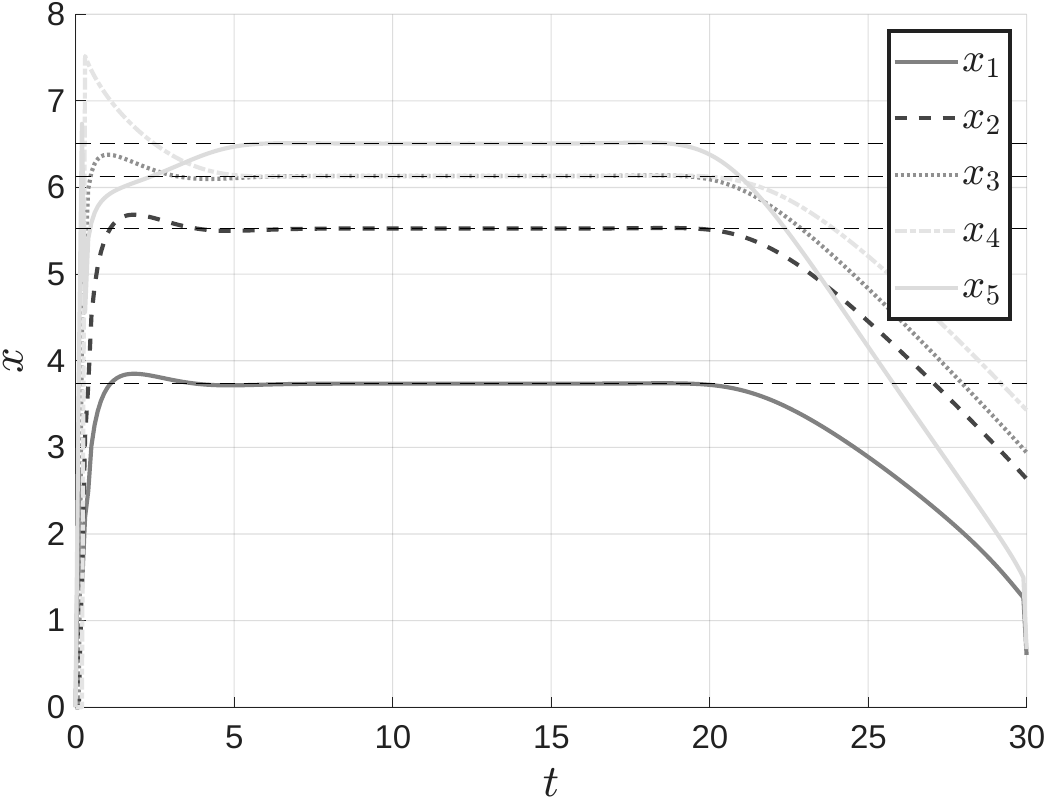}
			\caption[b]{$T = 30$}
		\end{subfigure}\hfil
		\begin{subfigure}[t]{0.45\linewidth}
			\centering
			\includegraphics[width=\linewidth]{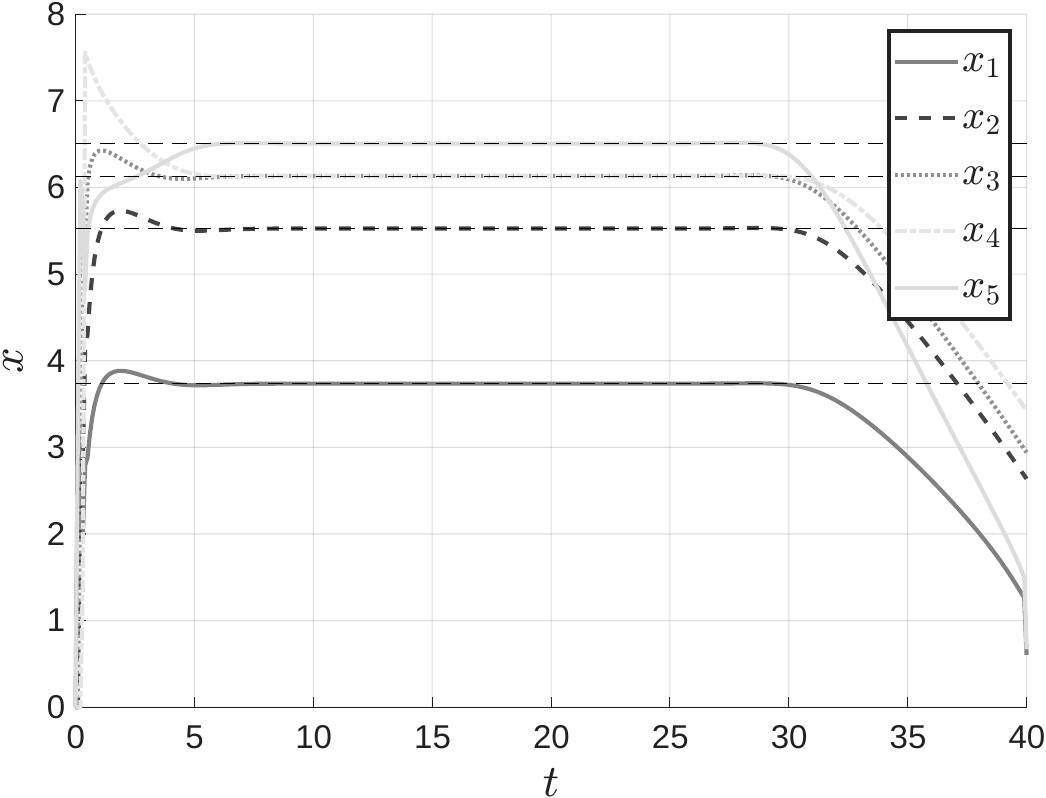}
			\caption[b]{$T = 40$}
			\label{fig:OCP_experiments_d}
		\end{subfigure}
		\caption{The state trajectories, $x^\star$, the solutions to $(\OCP_T^{\mathrm{heat}})$, with $x_0 = \zero$, $\bar u = 50$, $\Lambda_{12} = \Lambda_{23}=\Lambda_{34}=  \Lambda_{35}= 1$, $\alpha_1 = \alpha_2 = \alpha_3 = 1$ and different horizon lengths, $T$. The turnpike, $x^{\mathrm{tp}} \approx (3.74, 5.53, 6.13, 6.13, 6.51)^\top$, is indicated by horizontal dashed lines.}
	\end{figure}

\begin{figure}[htb] 
		\centering
		\begin{subfigure}[t]{0.45\linewidth}
			\centering
			\includegraphics[width=\linewidth]{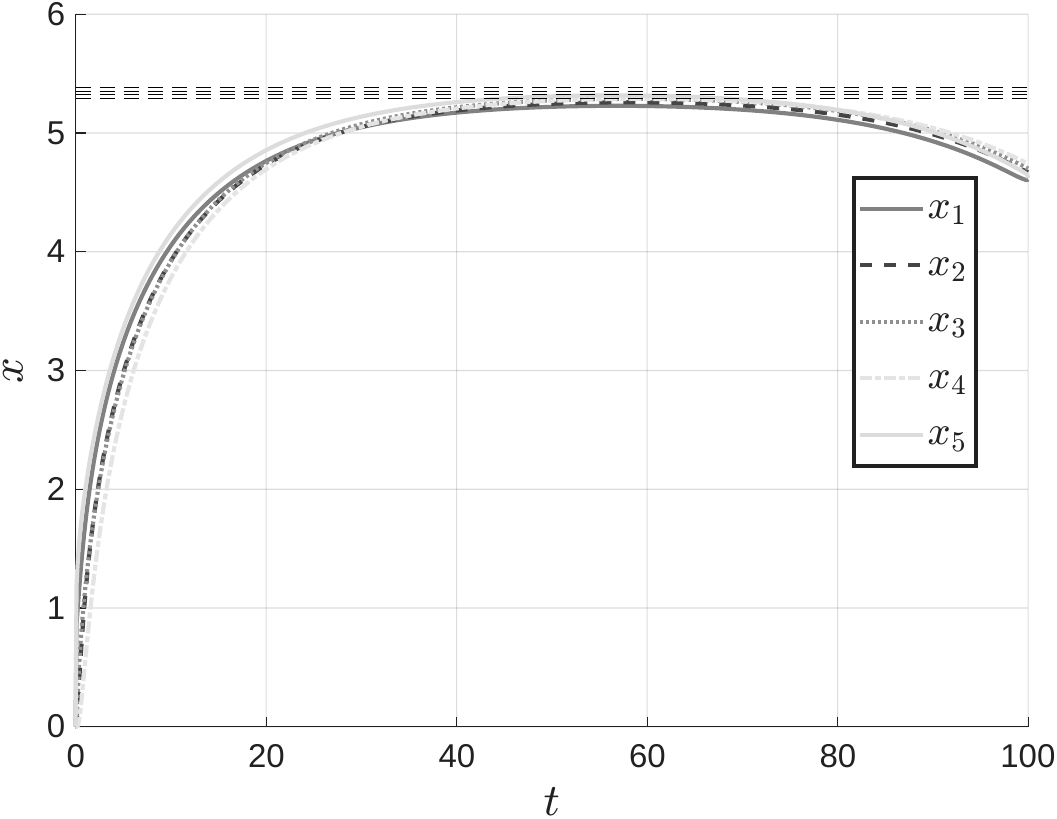}
			\caption{Parameters: $\Lambda_{12} = \Lambda_{34} = \Lambda_{35} = 0.1$, $\Lambda_{23} = 1$, $\alpha_1 = \alpha_3 = 1$, $\alpha_2 = 100$, $y^{\mathrm{ref}} = (1,5,10)^\top$, $\bar u = 50$, $T = 100$, $x^{\mathrm{tp}}~ \approx~(5.3, 5.32, 5.35, 5.35, 5.38)^\top$.}
			\label{fig:OCP_experiments_2_a}
		\end{subfigure}\hfil
		\begin{subfigure}[t]{0.45\linewidth}
			\centering
			\includegraphics[width=\linewidth]{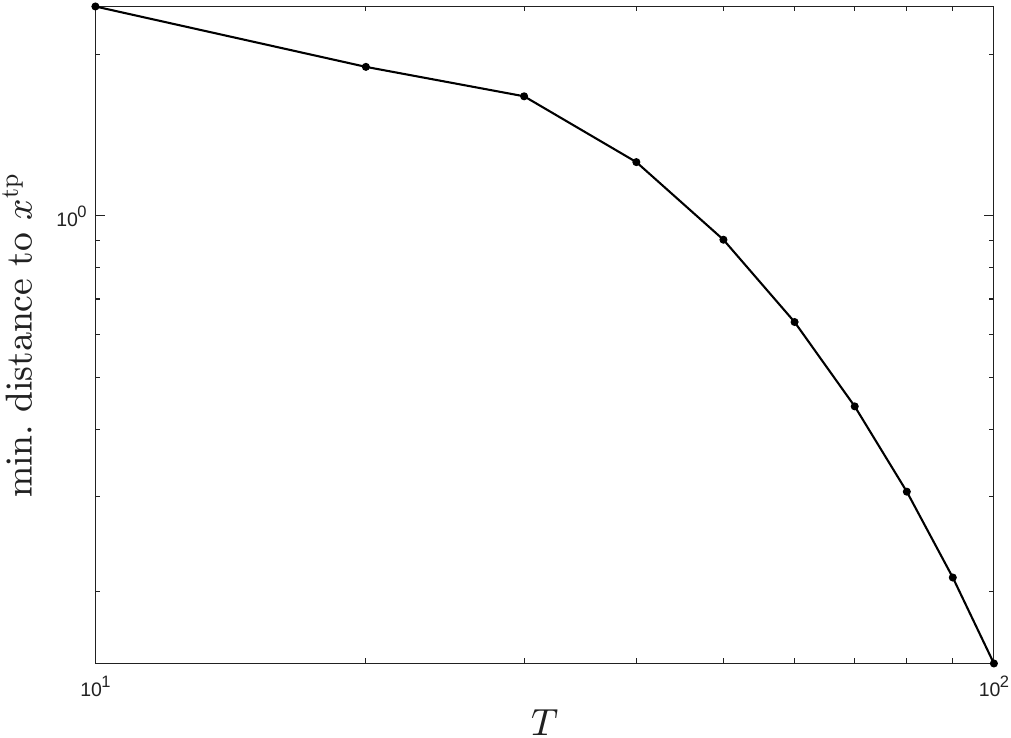}
			\caption[b]{Log-log plot of $\min_{t\in[0,T]} \| x^\star(t) - x^{\mathrm{tp}} \|$ versus $T$, parameters as in Figure~\ref{fig:OCP_experiments_2_a}. This suggests an exponential relationship.}
			\label{fig:OCP_experiments_2_b}
		\end{subfigure}
		\caption{The state trajectories, $x^\star$, the solution to $(\OCP_T^{\mathrm{heat}})$, showing the turnpike with large $\alpha_2$. Turnpike indicated with dashed horizontal lines.}
	\end{figure}

	Numerically we solve $(\OCP_T^{\mathrm{heat}})$ by discretising the dynamics via the direct Euler scheme with step size of $h=0.1s$ and solving the associated nonlinear program with the interior-point or SQP algorithm in Matlab's \texttt{fmincon} function. 
	Figures~\ref{fig:OCP_experiments_a}~--~\ref{fig:OCP_experiments_d} show the resulting optimal state trajectories with $x_0 = \zero$, $\bar u = 50$, $\Lambda_{12} = \Lambda_{23}=\Lambda_{34}=  \Lambda_{35}= 1$, $\alpha_1 = \alpha_2 = \alpha_3 = 1$ and different horizon  lengths. 
	These results strongly suggest that the solution to $(\OCP_T^{\mathrm{heat}})$ exhibits an \emph{integral state manifold turnpike}, which is defined as follows. 
	\begin{definition}[\cite{philipp2023optimal}]
		Consider $\OCPT$ with $\ell \in\C^1$, $\XM$ closed, and $f\in \C^1$. 
		We say that $\OCPT$ has the \emph{integral turnpike property} on $S_{\mathrm{tp}}\subset \mbbR^n$ with respect to a manifold $\mathcal{M}\subset\mbbR^n$ if for all compact $\mathbb{X}_0 \subset S_{\mathrm{tp}}$ there is a constant $C > 0$ such that for all $x_0\in\mathbb{X}_0$ and $T>0$, each optimal pair $(u^\star, x^\star)$ of $\OCPT$ with $x(0) = x_0$ satisfies,
		\[
		\int_0^T \dist^2(x^\star(t), \mathcal{M}) \, \dee t \leq C. 
		\]
	\end{definition}
	\noindent Thus, informally speaking, a turnpike is a set in the state space, close to which the optimal state $x^\star$ remains for the majority of the time.
	Here, we observe that the turnpike associated with $(\OCP_T^{\mathrm{heat}})$ is indeed only a point, $x^{\mathrm{tp}}\in\mbbR^5$, that, together with the control $u^{\mathrm{tp}}=(u^{\mathrm{tp}}_1, 0, 0, 0, u^{\mathrm{tp}}_5)^\top\in\mbbR^5$, solves the following nonlinear program,
	\begin{numcases}{(\NLP_{\mathrm{tp}})\quad}
	\min\limits_{(x,u)\in\mbbR^5 \times \mathbb{U}} & $\quad \mathcal{I}_{\mathrm{tp}}(x)$, \nonumber \\
	\mathrm{subject\  to:}	
	& $\quad f(x,u) = 0$,\nonumber
	\end{numcases}
	where $f$ is the right-hand side of the differential equation appearing in $(\OCP_T^{\mathrm{heat}})$ and
	\begin{align*}
	\mathcal{I}_{\mathrm{tp}}(x)  
	& = \biggl(\alpha_2\biggl[\sum\nolimits_{\substack{i\in[2:n] \\ j\in[i+1:n]}} \gamma_{(i,j)}(x) \left(\{S, H\}_{J^{(i,j)}}(x)\right)^2 \biggr]  + \alpha_3 \| C x- y^{\mathrm{ref}} \|^2\biggr)
	\end{align*}
	is the term appearing in the integral in \eqref{eq:integral_cost}, where we recall from Section~\ref{sec:app:heat_ex} that we use the indices $i$ and $j$ instead of $k$, for convenience.
	From the constraint $f(x^{\mathrm{tp}},u^{\mathrm{tp}}) = 0$ we see that $x^{\mathrm{tp}}$ is a controlled equilibrium point of the system. 
	We also observe that $\alpha_1$ does not affect the turnpike, only influencing how large the ``dip'' in the state at the end of the horizon is, as $\alpha_1$ penalises $H(x(T))$ in \eqref{eq:integral_cost}.  Recall that 
	\[
	\sum\nolimits_{\substack{i\in[2:n] \\ j\in[i+1:n]}} \gamma_{(i,j)}(x) \left(\{S, H\}_{J^{(i,j)}}(x)\right)^2 = \sum\nolimits_{\substack{i\in[2:n] \\ j\in[i+1:n]}} \Lambda_{ij}\tfrac{\left( e^{x_i} - e^{x_j} \right)^2}{e^{x_i}e^{x_j}}. 
	\]
	Therefore, if $\alpha_2$ is large with respect to $\alpha_3$ (i.e., entropy growth is heavily penalised) then the components of $x^{\mathrm{tp}}$ tend to be close to one another, whereas if $\alpha_3$ is large with respect to $\alpha_2 $ then $Cx^{\mathrm{tp}}$ tends to be closer to the references $y^{\mathrm{ref}}$.
	See Figure~\ref{fig:OCP_experiments_2_a}
	for the turnpike behaviour with large $\alpha_2$. 
	Finally, with the same parameters as those in Figure~\ref{fig:OCP_experiments_2_a}, Figure~\ref{fig:OCP_experiments_2_b} shows a log-log plot of the closest distance to the turnpike, $\min_{t\in[0,T]} \| x^\star(t) - x^{\mathrm{tp}} \|$, as a function of the horizon length $T$ suggesting that this dependence is exponential. 
	
	Although the paper \cite{philipp2023optimal} analysed the turnpike phenomena appearing in heat exchanger networks, the main result, \cite[Thm.~4.2]{philipp2023optimal}, establishes the existence of a turnpike under the quite strict condition that there is at least one \emph{thermodynamic equilibrium} in $C^{-1}(y^{\mathrm{ref}})$, that is, it requires that $C^{-1}(y^{\mathrm{ref}}) \cap \{x\in\mbbR^n :\sum\nolimits_{k=1}^N \gamma_k(x, H_x(x)) \left(\{S, H\}_{J_k}(x)\right)^2 = 0 \} \neq \emptyset.$
	For the heat exchanger network we considered in the numerical simulations of this section this assumption translates into $y^{\mathrm{ref}}_1 = y^{\mathrm{ref}}_4 = y^{\mathrm{ref}}_5$, a condition that is not satisfied in our examples. 
	Future research could build on the work conducted in \cite{philipp2023optimal} by trying to rigorously explain the occurrence of the observed turnpikes. 
	
	\subsection{Model predictive control of the heat exchanger network}
	
	We now apply model predictive control (MPC), see for example \cite{Rawlings2017,grune2017nonlinear}, to the heat exchanger network of the previous subsection. 
	MPC is an optimization-based feedback controller that iteratively solves a finite-horizon optimal control problem, see Algorithm~\ref{alg:mpc} for details.
	The state along the corresponding closed-loop trajectory is denoted by $x^{\mathrm{cl}}$. 
	\begin{algorithm}[h]
		\caption{MPC algorithm}
		\noindent \textbf{Input}: control horizon, $\delta >0$, prediction horizon, $T> 0$, initial state, $x^0 \in \mbbR^5$
		\begin{enumerate}
			\item Find the optimal control, $u^\star\in\U_T$, to $(\OCP_T^{\mathrm{heat}})$.
			\item Implement $u^\star(\boldsymbol \cdot)$, for $t \in [0, \delta]$.
			\item Set $x_0\leftarrow x(\delta;x_0,u^{\star})$ and go to step (1).
		\end{enumerate}
		\label{alg:mpc}
	\end{algorithm}
	
\begin{figure}[htb] 
		\centering
		\begin{subfigure}[t]{0.45\linewidth}
			\centering
			\includegraphics[width=\linewidth]{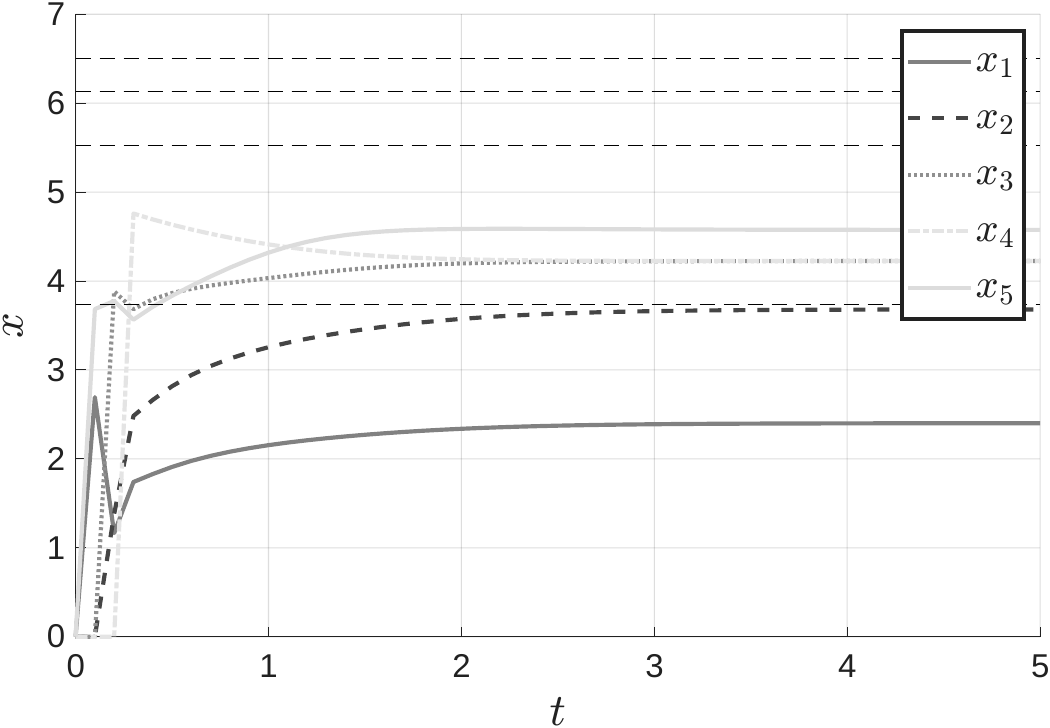}
			\caption{$T = 4$}
			\label{fig:MPC_experiments_a}
		\end{subfigure}\hfil
		\begin{subfigure}[t]{0.45\linewidth}
			\centering
			\includegraphics[width=\linewidth]{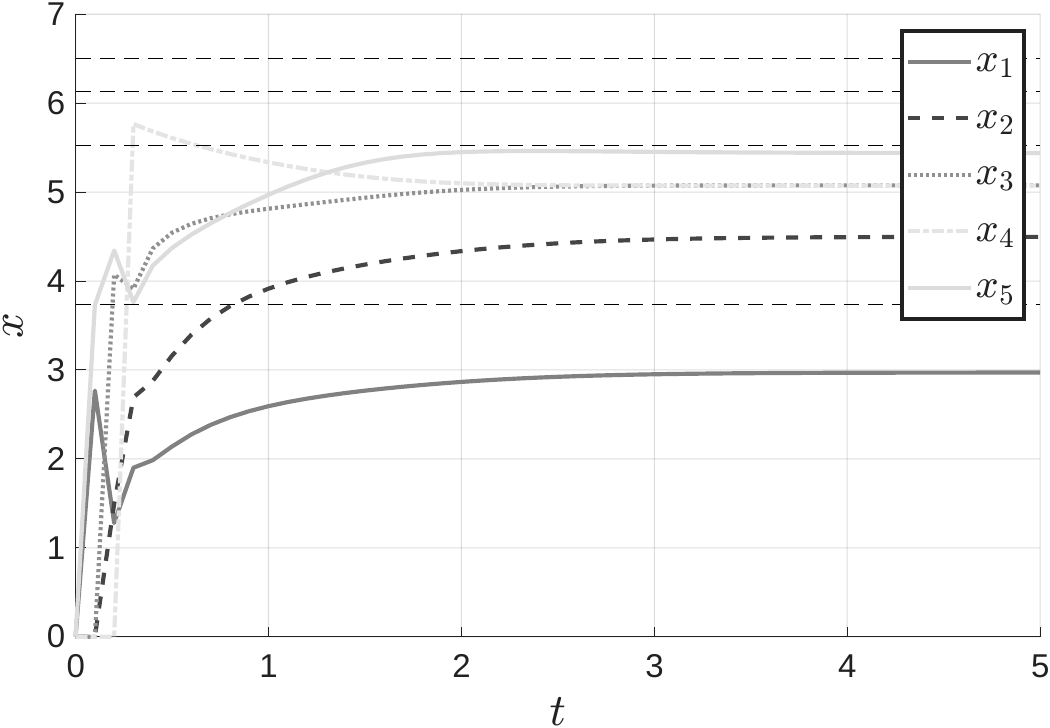}
			\caption[b]{$T = 6$}
		\end{subfigure}
		\begin{subfigure}[t]{0.45\linewidth}
			\centering
			\includegraphics[width=\linewidth]{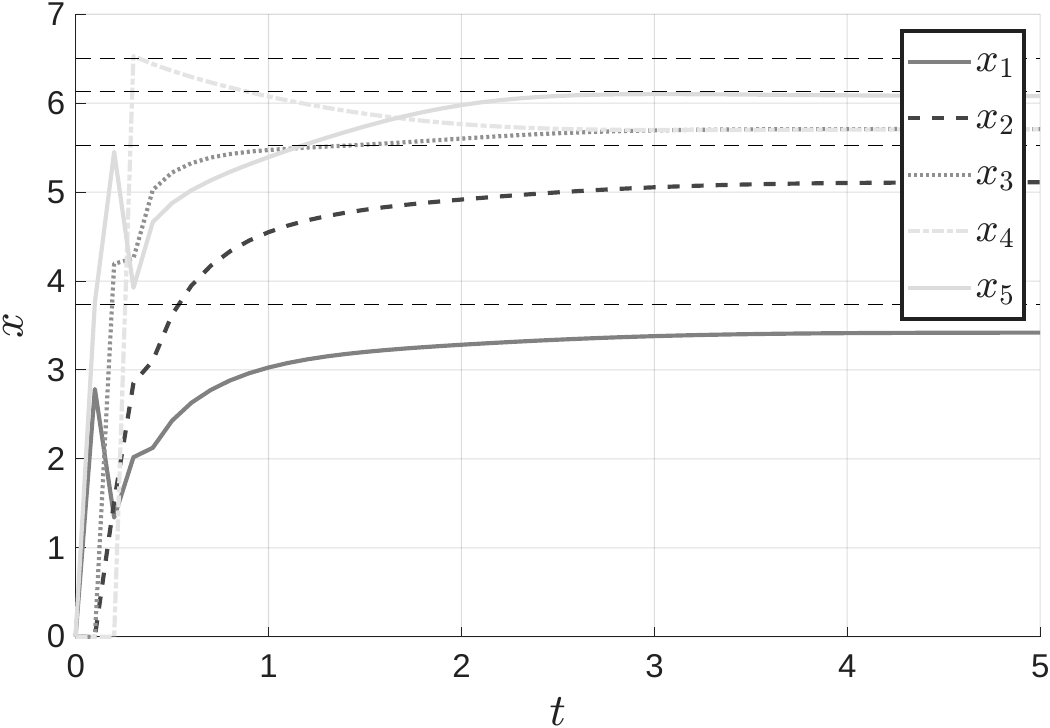}
			\caption[b]{$T = 8$}
		\end{subfigure}\hfil
		\begin{subfigure}[t]{0.45\linewidth}
			\centering
			\includegraphics[width=\linewidth]{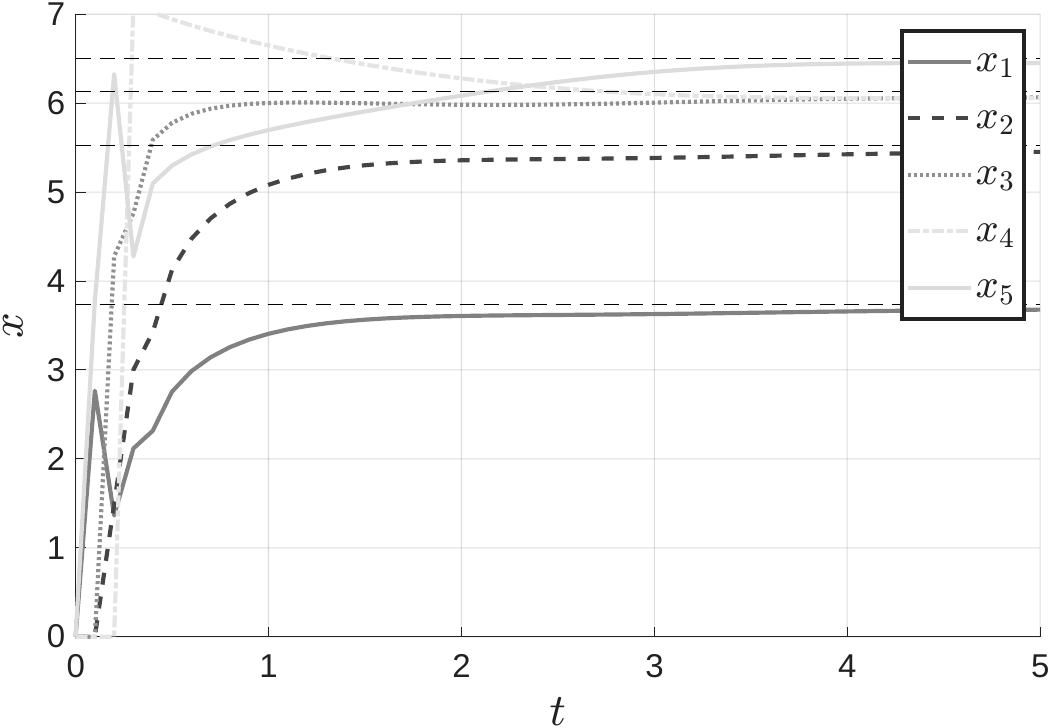}
			\caption[b]{$T = 10$}
			\label{fig:MPC_experiments_d}
		\end{subfigure}
		\caption{Closed-loop state trajectories, $x^{\mathrm{cl}}$, produced via MPC (as in Algorithm~\ref{alg:mpc}), with $x_0 =\zero$, $\bar u = 10$, $\Lambda_{12} = \Lambda_{23}=\Lambda_{34}=  \Lambda_{35}= 1$, $\alpha_1 = \alpha_2 = \alpha_3 = 1$, $\delta = 0.1s$, and different horizon lengths in $(\OCP_T^{\mathrm{heat}})$. 
			The closed-loop asymptotically settles at a steady-state, which approaches the turnpike, $x^{\mathrm{tp}} \approx (3.74, 5.53, 6.13, 6.13, 6.51)^\top$ (dashed horizontal lines), with increasing $T$.}
	\end{figure}

    \noindent Figures \ref{fig:MPC_experiments_a}~--~\ref{fig:MPC_experiments_d} show the closed loop $x^{\mathrm{cl}}$ obtained with different horizon lengths $T$ in $(\OCP_T^{\mathrm{heat}})$ and a control horizon $\delta$ = 0.1s (i.e., the same as the step size in the Euler scheme).	
    Instead of a dip at the end of the time interval, as in the OCP of the previous subsection, we observe that $x^{\mathrm{cl}}$ settles at a steady state, which is determined by $T$, and that this steady state approaches the turnpike as $T$ increases. 
    The reason for this lack of a dip is probably due to the control horizon, $\delta$, being very small (which is of course required for a well-designed feedback) resulting in this dip not having any influence on the long term behaviour of the closed-loop state.   
    	Future research could underpin these observations theoretically, drawing from established theory on the relationship between the MPC closed loop's steady state and turnpikes, see \cite{grune2016approximation}.
	
	\section{Conclusion and Outlook}\label{sec:conclusion}
	
	\noindent Exploiting the particular properties of coupled RI-PHSs, we have presented conditions under which the solution to an IVP remains bounded for any admissible control, from which existence of a global solution can be deduced. 
	We showed that this result implies the existence of a solution to a finite-horizon optimal control problem, if one additionally constrains the control to a compact set, the target is closed and if there exists an admissible state-control pair. 
	We verified that the conditions of the theoretical results are satisfied by networks of heat exchangers and a gas-piston system. Last, we provided a thorough case study for the network of heat exchangers in view of long-term behavior of solutions and MPC. 
    
	Future work could consider an intrinsic choice of the tracking term in \eqref{eq:special_cost} aiming at stabilization towards $x_{\mathrm{tp}}$. Replacing the tracking term $\|Cx-y_\mathrm{ref}\|$ by the availability function \eqref{eq:AvailabilityFunction} at $x^{\mathrm{tp}}$ leads to a cost functional only defined via thermodynamic properties of the systems, i.e., it does not necessitate the choice of an output matrix and a norm on the state space.

\bibliographystyle{plain}
\bibliography{bibliography}

\appendix
\section{Filippov's existence theorem}

\noindent For completeness we state Filippov's existence theorem for optimal control problems with integral cost (so-called Lagrange problems). 
We present a tailored variant of the theorem with time-invariant constraints and dynamics. 
The reader may consult \cite[Ch.9]{cesari1983optimization} for a more general presentation (with non-autonomous dynamics and constraints) and \cite[Ch.~4]{macky1982introduction} for an introductory treatment of existence theorems in optimal control. Consider the optimal control problem
\begin{align}
\min_{u\in\U_T, x\in\AC_T} & \quad\mathcal{I}(u,x) = \int_0^T \ell(x(s), u(s))\,\, \dee s, \nonumber\\
\mathrm{subject\,\, to:} & \quad \dot x (t) = f(x(t), u(t)),\quad \mathrm{a.e.}\quad t\in[0,T], \label{eq_Filipp_1}\\
& \quad x(t) \in \mathbb A, \quad \forall t\in[0,T], \label{eq_Filipp_2} \\
& \quad (x(0), x(T)) \in \mathbb B.  \label{eq_Filipp_3}
\end{align}
Here, $\AC_T$ denotes absolutely continuous functions mapping $[0,T]$ to $\mbbR^n$ and we let $\U_T=\Lcal([0,T], \mathbb U)$ being all Lebesgue-measurable functions mapping $[0,T]$ to $\mathbb U \subseteq \mbbR^m$. 

\begin{theorem}\cite[Thm.~9.3.i, Ch.9]{cesari1983optimization}\label{thm:filippov}
	Let $\mathbb{A}$ be compact, $\mathbb B$ be closed, $\mathbb M := \mathbb A \times \mathbb{U}$ be compact, and $f$ and $\ell$ be continuous on $\mathbb M$. 
	Assume that the set of admissible pairs, 
	\[
	\Omega := \{(x,u)\in\AC_T\times \U_T : \eqref{eq_Filipp_1}-\eqref{eq_Filipp_3}\,\, \mathrm{hold} \}
	\]
	is nonempty and that the set
	\[
	Z(x) := \{(z^0, z^\top)^\top\in\mbbR\times \mbbR^n : z = f(x,u), \,\, z^0 \geq \ell(x,u),\,\, u\in \mathbb{U} \}
	\]
	is a convex set in $\mbbR^{n + 1}$ for every $x \in \mathbb A$. 
	Then $\mathcal I$ has a global minimum in $\Omega$. 
\end{theorem}

\end{document}